\documentclass[11pt,twoside]%{cctart}
 {article}
\usepackage{mathrsfs}
\makeatother

\normalsize
\usepackage{amsfonts}
\usepackage[mathscr]{euscript}

\newcommand{\lbl}[1]{\label{#1}}

\newtheorem{theo}{Theorem}[section]
\newtheorem{prop}{Proposition}[section]
\newtheorem{lem}{Lemma}[section]
\newtheorem{remark}{Remark}[section]
\newtheorem{col}{Corollary}[section]
\newtheorem{defi}{Definition}[section]
\newcommand{\be}{\begin{equation}}
\newcommand{\ee}{\end{equation}}
\newcommand\bes{\begin{eqnarray}}
\newcommand\ees{\end{eqnarray}}
\newcommand{\bess}{\begin{eqnarray*}}
\newcommand{\eess}{\end{eqnarray*}}
\newcommand\bessd{\begin{eqnarray*}\left\{\begin{array}{ll}\smallskip}  \newcommand\eessd{\end{array}\right.\end{eqnarray*}}

\newcommand\ep{\varepsilon}
\newcommand\dd{\displaystyle}
\newcommand\vp{\varphi}

\newcommand\lk{\left}
\newcommand\rr{\right}

\newcommand\xxx{\bes\lk\{\begin{array}{ll}\smallskip}
\newcommand\zzz{\end{array}\right.\ees}

\def\theequation{\arabic{section}.\arabic{equation}}
\begin{document}
\setlength{\baselineskip}{16pt}

\begin{center}{\Large\bf  Initial and boundary blow-up problem}\\ [2mm]
{\Large\bf for $p$-Laplacian parabolic equation with general absorption} \footnote{The first author was supported by NSFC Grant 11371113; the second author was supported by
NUS AcRF Grant R-146-000-123-112; and the
third author was supported by Nantong Applied
Research Project Grant K2010042.}\\[4mm]
{\large Mingxin Wang$^{a,}$\footnote{Corresponding author. E-mail: mxwang@hit.edu.cn; Tel: 86-15145101503; Fax: 86-451-86402528}, \ \ Peter Y. H. Pang$^b$, \ \ Yujuan Chen$^c$} \\[2mm]
$^a$ {\small Natural Science Research Center, Harbin Institute of
  Technology, Harbin 150080, P. R. China.}\\[1mm]
$^b$ {\small Department of Mathematics, National University of Singapore, 10 Lower Kent Ridge Road,}\\
{\small Republic of Singapore 119076}\\[1mm]
$^c$ {\small Department of Mathematics, Nantong University, Nantong 226007, P. R. China}
 \end{center}

\begin{quote}
\noindent {\bf Abstract.} {\small In this article, we investigate the initial and boundary blow-up problem for the $p$-Laplacian parabolic equation $u_t-\Delta_p u=-b(x,t)f(u)$ over a smooth bounded domain $\Omega$ of $\mathbb{R}^N$ with $N\ge2$, where $\Delta_pu={\rm div}(|\nabla u|^{p-2}\nabla u)$ with $p>1$, and $f(u)$ is a function of regular variation at infinity. We study the existence and uniqueness of positive solutions, and their asymptotic behaviors near the parabolic boundary.}


\end{quote}

\section{Introduction and main results}
\setcounter{equation}{0} \thispagestyle{empty}
 {\setlength\arraycolsep{2pt}

Let $\Omega\subset\mathbb R^N\ (N\geq 2)$ be an open bounded domain with smooth boundary $\partial\Omega$, and $\Omega_T:=\Omega\times(0,T)$ with $0<T<\infty$. The aim of this paper is to study the $p$-Laplacian parabolic equation
 \bes
 u_t-\Delta_p u=-b(x,t)f(u),\ & (x,t)\in \Omega_T,\lbl{1.1}
 \ees
with blow-up initial and boundary values:
  \bes
 && u=\infty, \ \ \   (x,t)\in \partial\Omega\times(0,T),\lbl{1.2}\\
 && u=\infty, \   \ \ \ (x,t)\in\bar\Omega\times\{0\},
  \lbl{1.3}
 \ees
where $\Delta_pu={\rm div}(|\nabla u|^{p-2}\nabla u)$ with $p>1$,
$b(x, t)$ is a positive continuous function in $\Omega_T$ ($b(x,T)=0$ or $b(x,T)=\infty$ is allowed), and $f\in C^1([0,\infty))$ with $f(0)=0$ and $f'(u)>0$ for $u>0$.

Throughout this work, by (\ref{1.2})-(\ref{1.3}), we mean that
 \bess\lk\{\begin{array}{ll}\smallskip
 u(x,t)\to\infty\ \ & {\rm as}\ d(x)\to0\ {\rm uniformly \ for}\ t\in(0,T),\\
 u(x,t)\to\infty & {\rm as}\ t\to 0\ {\rm uniformly \ for }\ x\in\bar\Omega,
 \end{array}\rr.\eess
where, unless specified otherwise, $d(x)=d(x,\partial\Omega)={\rm
dist}(x,\partial\Omega)$ represents the distance from $x$ to
$\partial\Omega$ for $x\in\Omega$.

\begin{remark} The results of this paper remain valid if we add the term $a(x,t)u^{p-1}$ to the right hand side of the equation
$(\ref{1.1})$. For simplicity, we have not included this term.\end{remark}

We are interested in the
existence and uniqueness of positive weak solutions to
(\ref{1.1})--(\ref{1.3}), and the behavior of the solutions near the
parabolic boundary
 $$\Sigma_T:=\partial\Omega\times(0,T)\cup\bar\Omega\times\{0\}.$$

While there is an abundance of work -- going back to Bieberbach in 1916 -- on boundary blow-up for elliptic equations,
the corresponding investigation for parabolic equations has lagged behind. In 1994, Bandle et.~al.~\cite{BDD94}
studied the existence, uniqueness and asymptotic behavior near the
parabolic boundary of solutions to the autonomous parabolic boundary blow-up problem
  \bes\lk\{\begin{array}{ll}\smallskip
  u_t-\Delta\phi(u)=-f(u),\ \ &(x,t)\in\Omega\times(0,\infty), \\
  u=\infty,\ \ &(x,t)\in\partial\Omega\times(0,\infty)\cup\bar\Omega\times\{0\}.
  \end{array}\rr.\lbl{1.4}
  \ees
In particular, they proved that, under suitable
conditions on the functions $\phi$ and $f$,
  \bess&\dd\frac{u(x,t)}{w(t)}\to1\ \ \ {\rm as} \ (x,t)\to \Omega\times\{0\},&\\
 &\dd\frac{u(x,t)}{V(x)}\to1\ \ \ {\rm as} \ (x,t)\to \partial\Omega\times(0,\infty),&\eess
where $w(t)$ is a solution of
   \bes
   w'=-f(w),\ \ t>0; \ \ \ \ w(0)=\infty,\lbl{1.5}
   \ees
 and $V (x)$ is the unique solution to the elliptic boundary blow-up problem
  \bes
  \Delta\phi(v)=f(v), \ \ \ x\in \Omega;\ \ \
  v|_{\partial\Omega}=\infty.\lbl{1.6}
  \ees

In \cite{MV09},  Marcus and V\'{e}ron showed that if $f$ is super-additive, i.e.,
  \[f(u+v)\geq f(u)+f(v),\ \ \forall \ u,v\geq 0,\]
and satisfies
$$\int_{a}^{\infty}\frac{ds}{f(s)}<\infty, \qquad \int_1^{\infty}\frac{ds}{\sqrt{F(s)}}<\infty,$$
where $a$ is a non-negative constant such that $f(u)$ is positive and continuous when $u>a$ and ${\displaystyle F(s)=\int_0^sf(\tau)d\tau}$,
then there exists a maximal solution $\bar u(x,t)$ to (\ref{1.4}), and
  \bess
  \bar u(x,t)\leq w(t)+V(x),\ \ \ \bar u(x,t)\geq \max\{w(t),V(x)\},\ \ \forall\ (x,t)\in\Omega\times(0,\infty).
  \eess

For the non-autonomous case, very recently, motivated by a spatial-temporal degeneracy problem for the diffusive logistic equation used in population dynamics, Du et.~al.~\cite{DPP11} investigated the following problem:
   \bess\left\{\begin{array}{ll}\smallskip
 u_t-\Delta u=a(x,t)u-b(x,t)u^q, \ \ &  (x,t)\in\Omega\times(0,T),\\
 u=\infty, \  & (x,t)\in\partial\Omega\times(0,T)\cup\bar\Omega\times\{0\},
  \end{array}\right.
 \eess
where $q>1$, $a(x,t)$ and $b(x,t)$ are continuous functions in $\bar\Omega\times[0, T]$ and $\Omega\times[0, T]$, respectively, and $b(x,t)$ satisfies
  \[\alpha_1(t)d^\beta(x)\leq b(x,t)\leq \alpha_2(t)d^\beta(x),\ \ \ \forall\ (x,t)\in\Omega\times[0,T)\]
with $\beta>-2$, and $\alpha_1(t)$ and $\alpha_2(t)$ being positive continuous functions in
$[0, T)$. They also obtained existence, uniqueness and asymptotic behavior results.
Furthermore, under the extra condition that $b(x, t)\geq c(T-t)^\theta
d^\beta(x)$ for some constants $c > 0, \theta> 0$ and $\beta>-2$,
they showed that the positive solution that exists stays bounded in any compact subset
of $\Omega$ as $t$ increases to $T$, and hence solves the equation
up to $t = T$.

Related problems have also been studied by \cite{AV1, AV2, BIV10, LP11} and \cite{SP12}. Especially, in \cite{LP11} the authors proved the existence of large solutions for the problems
 \bessd
 u_t-{\rm div}\,a(x,t,u,\nabla u)+g(x,t,u,\nabla u)=f(x,t), \ \ &  (x,t)\in\Omega\times(0,T),\\ \smallskip
 u=u_0,\ &  (x,t)\in\Omega\times\{0\},\\
 u=\infty,&(x,t)\in\partial\Omega\times(0,T),
 \eessd
where ${\rm div}\,a(x,t,u,\nabla u)\approx\Delta_pu$, $g(x,t,u,\nabla u)\approx u|\nabla u|^q$ with $p-1<q\leq p$, and $u_0\in L^1_{\rm loc}(\Omega)$, $f\in L^1(0,T;\,L^1_{\rm loc}(\Omega))$ with $f^-\in L^1((0,T)\times\Omega)$. In \cite{SP12}, the existence and uniqueness of entropy large solutions was discussed for the following problem
 \bessd
 u_t-\Delta_pu=0, \ \ &  (x,t)\in\Omega\times(0,T),\\ \smallskip
 u=u_0,\ &  (x,t)\in\Omega\times\{0\},\\
 u=\infty,&(x,t)\in\partial\Omega\times(0,T),
 \eessd
where $1\leq p<2$, $u_0\in L^1_{\rm loc}(\Omega)$  ($u_0\in L^1(\Omega)$ if $p=1$) is a nonnegative function.

\vskip 4pt Motivated by the above works, in this paper, we study the problem (\ref{1.1})--(\ref{1.3}). We are able to extend some of the results of \cite{BDD94,DPP11,MV09}. Our method refers to Karamata's regular variation theory \cite{BGT87}, which has been used by many authors in elliptic boundary blow-up problems.

We briefly recall some key notions of Karamata's theory; more can be
found in the Appendix.

 A measurable function $R: [A,\infty)\to
(0,\infty)$, for some $A>0$, is called {\em regularly varying} at
infinity of index $\rho\in\mathbb R$, for short $R\in RV_\rho$, if
$\dd\lim_{u\to\infty}\frac{R(\xi u)}{R(u)}=\xi^\rho,
   \ \ \forall \, \xi>0.$
When the index $\rho$ is zero, we call the function $R$ {\em slowly
varying} at infinity.

Following {\rm\cite{CR02}} (see also \cite{Mo07}), we denote by
$\mathcal {K}_\ell$ the set of all positive, monotonic functions $k\in
C^1(0, \mu)\cap L^1(0, \mu)$ that satisfy
 \[\lim_{s\to 0^+}\left(\frac{K(s)}{k(s)}\right)'=\ell\in(0,\infty),\]
where $K(s)=\int_0^sk(\theta){\rm d}\theta$ and $\mu\geq{\rm
diam}(\Omega)$. For any $k \in \mathcal{K}_\ell$, it is clear that $\dd\lim_{s\to 0^+}\frac{K(s)}{k(s)}= 0$ and $\dd\lim_{s\to 0^+}\frac{K(s)k'(s)}{k^2(s)}=1-\ell$.
Moreover, $0<\ell\leq1$ if
$k$ is non-decreasing, and $\ell\geq1$ if $k$ is non-increasing.

\vskip 4pt
With regard to (\ref{1.1}), we shall often make the following assumptions:

\vskip 2pt

$(F_1)$ ~~$f\in RV_\rho$ with $\rho>p-1$;

\vskip 2pt

$(F_2)$ ~~The function $s\mapsto s^{-(p-1)}f(s)$ is increasing in
 $(0,\infty)$;

\vskip 2pt

$(B)$ ~~There exist a function $k\in \mathcal {K}_\ell$ and
two positive continuous functions $\alpha_1(t)$ and $\alpha_2(t)$
defined on $[0, T)$, such that
  \bess
 \alpha_1(t)k^p(d(x))\leq b(x,t)\leq \alpha_2(t)k^p(d(x)),\ \ \forall \ (x,t)\in\Omega_T,
  \eess
where $\alpha_1(T)=0$ or $\alpha_2(T)=\infty$ may occur.

\begin{remark}\lbl{r2.z1} If we assume that both $\alpha_1(t)$ and $\alpha_2(t)$ are positive and continuous on $[0,T]$, then one can replace $\Omega_T$ by $Q_T:=\Omega\times(0,T]$ and the problem can be discussed in $Q_T$. Moreover, Theorem $\ref{t1.1}$ below will then also hold true for $t^*=T$.
\end{remark}

For notation, let $\phi$ be the function defined uniquely by
  \bes
  \int_{\phi(t)}^\infty\frac{ds}{(p'F(s))^{1/p}}
  =t,\ \ \ t>0,
 \lbl{1.10}\ees
where $F(t)=\int_0^t f(s){\rm d}s$ and $p'=\frac{p}{p-1}$. It is easily seen that $\phi(0)=\infty$. Further, let $\xi(t)$ be the unique positive solution of $(\ref{1.5})$ and $\xi^*(t)$ be the unique positive solution of
 \bes
 (\xi^*)'=-f^*(\xi^*), \ \ \ t>0; \ \ \ \xi^*(0)=\infty,
 \lbl{1.z1}\ees
with $f^*(s)=\lk(k\circ K^{-1}\circ\phi^{-1}(s)\rr)^pf(s)$.

\begin{theo}\lbl{t1.1} ~ {\rm(i)} \ Let the conditions $(F_1)$, $(F_2)$
and $(B)$ hold. Suppose that
  $$\rho>\max\left\{1, \,p-1,\, p-1-(p-2)/\ell\right\}.$$
Then the problem
$(\ref{1.1})$-$(\ref{1.3})$ has a maximal positive solution
$\overline{u}$ and a minimal positive solution $\underline{u}$, in
the sense that any positive solution $u$ of $(\ref{1.1})$-$(\ref{1.3})$
satisfies $\underline{u}\leq u\leq\overline{u}$. Moreover, the minimal positive solution $\underline{u}$ is non-increasing in $t$. Furthermore, for any given $t^*\in(0,T)$, there is a constant $C>0$, depending on $t^*$, such that the maximal positive solution $\overline{u}$ satisfies
   \bes
 \overline{u}(x,t)\leq\left\{\begin{array}{ll}\smallskip C\left[\xi(t)+\phi(K(d(x)))\right], & {\rm if}\ k\ \mbox{\rm is\ non-increasing},\\[1mm]
  C\left[\xi^*(t)+\phi(K(d(x)))\right], \ \ &{\rm  if}\
   k\ \mbox{\rm is\ non-decreasing},
   \end{array}\right.\lbl{1.11}
   \ees
for all $(x,t)\in\Omega\times(0,t^*]$.

{\rm(ii)} \ Assume that in addition $f$ satisfies the following condition:
 \vspace{-2mm}\begin{quote} $(C)$ \ There is a constant $l>\max\{1,p-1\}$ such that $f(u)\geq\varepsilon^{-l}f(\varepsilon u)$ for all $u>0$ and $0<\varepsilon\ll 1$.
  \vspace{-2mm}\end{quote}
Then, for any given $t^*\in(0,T)$, there is a constant $c>0$, depending on $t^*$, such that the minimal positive solution $\underline{u}$ satisfies
  \bes \underline{u}(x,t)\geq
 \left\{\begin{array}{ll}\smallskip c\left[\xi^*(t)+\phi(K(d(x)))\right],\ \ & {\rm if}\ k\ \mbox{\rm is\ non-increasing},\\[1mm]
  c\left[\xi(t)+\phi(K(d(x)))\right], \ \ &{\rm  if}\ k\ \mbox{\rm is\ non-decreasing},\end{array}\right.
  \lbl{1.12} \ees
for all $(x,t)\in\Omega\times(0,t^*]$.
  \end{theo}

\begin{remark}\lbl{rr1.2} \ {\rm(i)} \ If there is a constant $l>\max\{1,p-1\}$ such that the function $f(u)/u^{l}$ is increasing for $u>0$, then the condition $(C)$ holds.

{\rm(ii)} \  Set $q=\rho-(\rho-p+1)(1-\ell)$. Under the conditions of Theorem $\ref{t1.1}$, we have $f^*\in RV_{q}$ $($see $(\ref{b.1})$ below$)$.

{\rm(iii)} \ Clearly, $\rho>q$ if \ $0< \ell< 1$; \ $\rho= q$ if \ $\ell= 1$; \ $\rho <q$ if \ $\ell>1$.
As $\rho>\max\big\{1,\, p-1, \,p-1-(p-2)/\ell\big\}$, we have
$q>\max\{1,p-1\}$.
\end{remark}

To simplify notation, we denote
 \[r=\frac{\rho+1}{\rho+1-p}.\]

\begin{theo}\lbl{t1.2} ~Under the assumptions of Theorem $\ref{t1.1}$, let $u(x, t)$ be any positive solution of $(\ref{1.1})$-$(\ref{1.3})$. Then the following hold:

{\rm(i)}  For any fixed $t_0\in (0, T)$ and $y\in\partial\Omega$, we have
 $$\dd\lim_{\Omega\ni x\to y}\frac{u(x,t_0)}{\phi (K(d(x)))}
 =\left(\frac{r+\ell-1}{r\beta(y,t_0)}\right)^{\frac{r-1}{p}}$$
provided that $\beta(x,t):=\frac{b(x,t)}{k^p(d(x))}$
can be extended to a continuous function on $\bar\Omega\times
(0, T)$.

{\rm (ii)} For any fixed $x_0\in\Omega$, let $\tau(t)$ be the unique positive solution of
 \bes
   \tau'=-b(x_0,0)f(\tau),\ \ t>0; \ \ \ \tau(0)=\infty.\lbl{1.13}
   \ees
Then
  \bes
  \dd\limsup_{t\to0}\frac{u(x_0,t)}{\tau(t)}\leq 1.
  \lbl{1.15}\ees
If in addition $p>2N/(N+2)$ and $f(s)/s$ is increasing for $s>0$, then
 \bes
 \dd\liminf_{t\to0}\frac{u(x_0,t)}{\tau(t)}\geq 1.
 \lbl{1.16}\ees
\end{theo}

%\begin{theo}\lbl{t1.3}
%~??{\bf Under the assumptions of Theorem {\rm\ref{t1.1}}, if $p\geq 2$ and the function %$k(s)\equiv 1$ in the condition $(B)$, then the problem {\rm(\ref{t1.1})-(\ref{t1.3})} has a %unique positive solution.}
%\end{theo}

\begin{theo}\lbl{t1.3}
~Under the assumptions of Theorem {\rm\ref{t1.1}}, if $p=2$ and $k(s)=1$,
and $f(u)$ is convex in $(0,\infty)$, then {\rm(\ref{t1.1})-(\ref{t1.3})} has a unique positive solution.
\end{theo}

This paper is organized as follows. In Section 2, we
prove the comparison principle. Section 3 is devoted to prove Theorem \ref{t1.1}. The proofs of asymptotic behavior and uniqueness (Theorems \ref{t1.2} and \ref{t1.3}) will be given in Section 4. The last section (Appendix) contains three parts: (i)
state and prove some relevant results of the Karamata's regular variation theory which will be used in the text (not all of which are readily available in the literature). Especially, Lemmas \ref{l2.5}-\ref{l2.8} play an important role in the proofs of Theorem \ref{t1.1} and Theorem \ref{t1.2}(i); (ii) prove some results on the unique solution of (\ref{1.5}), which will be used in the proof of Theorem \ref{t1.2}(ii); (iii) state some results on the corresponding elliptic boundary blow-up problem.

\section{Preliminaries}
\setcounter{equation}{0}

The main aim of this section is to prove the key comparison principle that is crucial to this paper. While the comparison principle is, in a sense, known, we believe a careful proof is useful to clarify the different versions that appear in the literature.

%This paper also makes substantial use of the Karamata regular variation theory. %However, for the sake of a more streamlined flow of the paper, instead of including %the relevant results here (not all of which are readily available in the literature), %we have relegated them to the Appendix.

We first establish a notation: If $\varphi\in C^\infty(\Omega_T)$ and ${\rm supp}\,\vp\subset\subset \Omega_T$, i.e. $\vp$ is zero near the parabolic boundary $\Sigma_T:=\partial\Omega\times(0,T)\cup\bar\Omega\times\{0\}$ of $\Omega_T$, we write
$\vp\in C_\bullet^\infty(\Omega_T)$.

\begin{defi}\lbl{dd2.1} \ A weak lower $($upper$)$ solution of the equation $(\ref{1.1})$ is a measurable function $u(x, t)$ such that
$$u\in C(t_0,T;\, L^2(\Omega'))\cap L^p(t_0,T;\, W^{1,p}(\Omega'))\cap L^\infty(\Omega'\times(t_0,T)),\ \ u_t\in L^2(\Omega'\times(t_0,T))$$
for any $0<t_0<T$ and any compact subset $\Omega'$ of $\Omega$;
and
\bes
  \int_{\Omega} u\varphi{\rm d}x+\int_{0}^t\int_{\Omega}\{|\nabla u|^{p-2}\nabla u\cdot\nabla\varphi+bf(u)\varphi\}{\rm d}x{\rm d}\tau
  \leq(\geq)\dd\int_{0}^t\int_{\Omega}u\varphi_t{\rm d}x
  {\rm d}\tau, \ \ \forall \ 0<t<T.
 \lbl{2.1}
  \ees
for all test function $\varphi\in C_\bullet^\infty(\Omega_T)$, $\vp\geq 0$ in $\Omega_T$.

A function $u$ that is both a lower solution
and a upper solution is a weak solution of the equation {\rm(\ref{1.1})}.
\end{defi}

\begin{prop}\lbl{p2.1}{\rm(Comparison Principle)}
Let $f\in C[0,\infty)$ be a non-negative function, and $
b(x, t)\in C(\Omega_T)$ be a non-negative and non-trivial function. Assume that $ u_1, u_2\in
C^1(\Omega_T)$ are weak upper and lower solutions of equation
$(\ref{1.1})$ respectively, that are positive in $\Omega_T$. If $f(s)$ is non-decreasing for $s\in \big(\inf_{\Omega_T}\{u_1,u_2\},\,\sup_{\Omega_T}\{u_1,u_2\}\big)$, and $u_1, u_2$ satisfy
  \bes
  \limsup_{(x,t)\to\Sigma_T}(u_2-u_1)\leq 0,
  \lbl{2.3c}\ees
then $u_1\geq u_2$ in $\Omega_T$.
\end{prop}

{\bf Proof}. ~The proof refers to the corresponding elliptic case (\cite{DG03}), and makes use of \cite[Lemma 2.1]{D98}. Let $\varphi\in C_\bullet^\infty(\Omega_T)$ be a non-negative function. Then we have
 \bes
 &&\dd\int_{\Omega}(u_2-u_1)\vp{\rm d}x+
 \int_{0}^t\int_{\Omega}\big(|\nabla u_2|^{p-2}\nabla u_2-|\nabla u_1|^{p-2}\nabla u_1\big)\cdot\nabla \varphi{\rm d}x{\rm d}\tau\nonumber\\
 &+&\int_{0}^t\int_{\Omega}b[f(u_2)-f(u_1)]\vp{\rm d}x
  {\rm d}\tau
 \leq\int_{0}^t\int_{\Omega}(u_2-u_1)\vp_t{\rm d}x{\rm d}\tau , \ \ \forall \ 0<t<T.
  \lbl{2.4}\ees
For any $0<\varepsilon<1$, let $v=[u_2-(u_1+\varepsilon)]_+$
where $u_+:=\max\{u,0\}$. By the assumption (\ref{2.3c}),
 \[\limsup_{(x,t)\to\Sigma_T}(u_2-u_1-\varepsilon)\leq \limsup_{(x,t)\to\Sigma_T}(u_2-u_1)-\varepsilon\leq-\varepsilon.\]
We can choose $t_\varepsilon\in(0,T)$ and $\Omega(\varepsilon)\subset\subset\Omega$ with $t_\varepsilon\to 0$ and $\Omega(\varepsilon)\to\Omega$ as $\varepsilon\to 0$, such that
$v=0$ in $\Omega_T\setminus\Omega(\varepsilon)\times(2t_\varepsilon, T)$ and
  \[v\in W^{1,2}(t_\varepsilon,T;\, L^2(\Omega(\varepsilon)))\cap L^p(t_\varepsilon,T;\, W^{1,p} (\Omega(\varepsilon)))\cap L^\infty(\Omega(\varepsilon)\times(t_\varepsilon,T)).\]
It follows that $v$ can be
approximated arbitrarily closely in the norm of
 $$W^{1,2}(t_\varepsilon,T;\, L^2(\Omega(\varepsilon)))\cap L^p(t_\varepsilon,T;\, W^{1,p} (\Omega(\varepsilon)))\cap L^\infty(\Omega(\varepsilon)\times(t_\varepsilon,T))$$
by $C_\bullet^\infty(\Omega_T)$ functions. Thus (\ref{2.4}) holds with $\varphi$ replaced by $v$. For any given $t_\varepsilon<s<T$, denote
 $$D^\varepsilon_s=\big\{(x,t)\in\Omega(\varepsilon)\times(t_\varepsilon, s]:\,u_2(x,t)>u_1(x,t)+\varepsilon\big\}, \ \ C^\varepsilon_s=\{(x,t)\in D^\varepsilon_s: t= s\}.$$
To simplify the notation we write $w=u_2-u_1$. Then for any fixed $t_\varepsilon<t<T$, we have
  \bes
 &&\dd\int_{C^\varepsilon_t}w(w-\varepsilon)_+{\rm d}x+\int_{D^\varepsilon_t}\big(|\nabla u_2|^{p-2}\nabla u_2-|\nabla u_1|^{p-2}\nabla u_1\big)\cdot\nabla(u_2-u_1){\rm d}x\nonumber\\
 &&+ \int_{D^\varepsilon_t}b\left(f(u_2)-f(u_1)\right)(w-\varepsilon)_+{\rm
 d}x{\rm d}\tau
 \leq\int_{D^\varepsilon_t}w[(w-\varepsilon)_+]_t{\rm d}x{\rm d}\tau.
  \lbl{2.5}\ees
It is obviously that the third term in the left hand side of (\ref{2.5}) is non-negative since $u_2>u_1$ in $D^\varepsilon_t$ and $f$ is non-decreasing and $b$ is positive. By \cite[Lemma 2.1]{D98}, we see that the second term in the left hand side of (\ref{2.5}) is also non-negative. Therefore
 \bes
 \dd\int_{C^\varepsilon_t}w(w-\varepsilon)_+{\rm d}x
 \leq\int_{D^\varepsilon_t}w[(w-\varepsilon)_+]_t{\rm d}x{\rm d}\tau.
  \lbl{2.5a}\ees
Noting that $C^\varepsilon_t\subset\Omega(\varepsilon)\times\{t\}$ and $w(x,t)\leq \varepsilon$ in $(\Omega(\varepsilon)\times\{t\})\setminus C^\varepsilon_t$, and $D^\varepsilon_t\subset
\Omega(\varepsilon)\times(t_\varepsilon, t]$ and $w(x,\tau)\leq\varepsilon$ in $\Omega(\varepsilon)\times(t_\varepsilon, t]\setminus D^\varepsilon_t$, we have
 \bess
 \int_{C^\varepsilon_t}w(w-\varepsilon)_+{\rm d}x&=&\int_{C^\varepsilon_t}(w-\varepsilon)_+^2{\rm d}x
 +\varepsilon\int_{C^\varepsilon_t}(w-\varepsilon)_+{\rm d}x\\
  &=&\int_{\Omega(\varepsilon)}(w-\varepsilon)_+^2{\rm d}x
  +\varepsilon\int_{\Omega(\varepsilon)}(w-\varepsilon)_+{\rm d}x,\\
 \int_{D^\varepsilon_t}w[(w-\varepsilon)_+]_t{\rm d}x{\rm d}\tau
  &=&\frac 12\int_{t_\varepsilon}^t\int_{\Omega(\varepsilon)}[(w-\varepsilon)_+^2]_t{\rm d}x{\rm d}\tau+\varepsilon\int_{t_\varepsilon}^t\int_{\Omega(\varepsilon)}[(w-\varepsilon)_+]_t{\rm d}x{\rm d}\tau\\
  &=&\frac 12\int_{\Omega(\varepsilon)}(w-\varepsilon)_+^2{\rm d}x+\varepsilon\int_{\Omega(\varepsilon)}(w-\varepsilon)_+{\rm d}x.
  \eess
This combined with (\ref{2.5a}) yields
 \[\int_{\Omega(\varepsilon)}(w(x,t)-\varepsilon)_+^2{\rm d}x=0, \ \ \ \forall \ t_\varepsilon<t<T,\]
which implies that $w(x,t)\leq\varepsilon$, i.e., $u_2(x,t)\leq u_1(x,t)+\varepsilon$ in $\Omega(\varepsilon)\times(t_\varepsilon, T]$.
Noting that $t_\varepsilon\to 0$ and $\Omega(\varepsilon)\to\Omega$ as $\varepsilon\to 0^+$, we conclude $u_2(x,t)\leq u_1(x,t)$ a.e. in $\Omega_T$ and complete the proof.\ \ \ \fbox{}

\vskip 8pt

For ease of reference, we end this section by recalling the following comparison principle for the corresponding elliptic problem, which can be derived from the characterizations of the maximum principle in \cite{GS98} and Proposition 2.2 in \cite{DG03, LPW12}.

\begin{prop}\lbl{p2.2}{\rm(Comparison Principle)}
~Suppose that $D$  is a bounded domain in $\mathbb R^N$, and $\beta(x)$ is a continuous function in $D$ with $\beta(x)\geq 0, \beta(x)\not\equiv 0$. Let
$u_1,\,u_2\in C^1(D)$ be positive in $D$ and satisfy in the sense of
distribution
 \bess
 -\Delta_p u_1+\beta(x)g(u_1)\geq 0
 \geq -\Delta_p u_2+\beta(x)g(u_2)
 \eess
and
 \[\limsup_{d(x,\,\partial D)\to 0}(u_2-u_1)\leq 0,\]
where $g\in C([0,\infty),[0,\infty))$. If furthermore we assume
that $g(s)/s^{p-1}$ is increasing for $s\in
\big(\inf_D\{u_1,u_2\},\,\sup_D\{u_1,u_2\}\big)$, then $u_1\geq u_2$
in $D$.
\end{prop}

\section{Maximal and minimal positive solutions}
\setcounter{equation}{0}

In this section, we give the proof of Theorem \ref{t1.1}. We shall divide the proof into five steps. Some of the techniques used are based on those found in \cite{DPP11} and \cite{CW11}, though the adaptation to our setting is not straightforward.

\vskip 4pt{\bf Step 1: Construction of upper solution.}

This is the key step in the proof. By Theorem \ref{t2.2}(i), the problem
 \xxx
  \Delta_p z= k^p(d(x))f(z), \ \ &x\in \Omega,\\
  z=\infty,\ \ &x\in \partial\Omega
 \lbl{3.1}\zzz
has a unique positive solution $z(x)$, and there are positive constants $c_1$ and $c_2$ such that
 \bes
 c_1\phi(K(d(x)))\leq z(x)\leq c_2 \phi(K(d(x))),\ \ \ x\in\Omega.\lbl{3.2}
 \ees

Case (1): {\it $k$ is non-increasing.}

For arbitrarily small $\varepsilon >0$, since $\alpha_1(t)>0$ in
$[0, T-\varepsilon]$, we may assume that $\alpha_1(t) \geq
\alpha_\varepsilon$ on $[0, T-\varepsilon]$ for some constant
$\alpha_\varepsilon>0$.
Let $\xi(t)$ be the unique positive solution of (\ref{1.5}). By the assumption on $k$, we can find $c>0$ such that $c k^p(d(x))\geq1$ in $\Omega$. It follows that
$\xi'(t)\geq-c k^p(d(x))f(\xi(t))$ in $\Omega_T$.
By Lemma \ref{l2.5}, we can find $\Lambda>1$ sufficiently large such that
$(c\Lambda+\Lambda^{p-1})f(\xi+z)<\alpha_\varepsilon f(\Lambda\xi+\Lambda z)$.
Since $f$ is increasing, the function $\overline{\mathbf{u}}(x,t)=\Lambda[\xi(t)+z(x)]$ satisfies (here, $d=d(x)$)
  \bes
  \overline{\mathbf{u}}_t-\Delta_p \overline{\mathbf{u}} &\geq&-c\Lambda k^p(d)f(\xi)-\Lambda^{p-1} k^p(d)f(z)\nonumber\\
 &\geq&-c\Lambda k^p(d) f(\xi+z)-\Lambda^{p-1} k^p(d) f(\xi+z)\nonumber\\
 &\geq&-\alpha_\varepsilon k^p(d)f(\overline{\mathbf{u}})\nonumber\\
 &\geq&-b(x,t)f(\overline{\mathbf{u}}), \ \ \ \ (x,t)\in \Omega\times(0,T-\varepsilon].
 \lbl{3.z1}\ees

 Case (2):  {\it $k$ is non-decreasing.}

Let $c_1$ and $c_2$ be given by (\ref{3.2}). Noting that $\phi$ is decreasing, $K$ is increasing and $k$ is non-decreasing, it follows that
 \bess
 \left\{\begin{array}{l} k(d)\geq k\circ K^{-1}\circ\phi^{-1}(c^{-1}_1s)\ \ \ \mbox{if} \ \ \phi(K(d))\leq c_1^{-1} s,\\[2mm]
 k(d)\leq k\circ K^{-1}\circ\phi^{-1}(c^{-1}_2s)\ \ \ \mbox{if} \ \ \phi(K(d))\geq c_2^{-1}s.\end{array}\right.\eess
Set $f_i(s)=\lk(k\circ K^{-1}\circ\phi^{-1}(c^{-1}_is)\rr)^pf(s)$. Then $f_i(s)\in RV_{q}$ by (\ref{b.1}), and
 \bes
 f_1(s)\leq k^p(d)f(s) \ \ \mbox{when} \ s\geq c_1\phi(K(d)), \ \ \ \
 f_2(s)\geq k^p(d)f(s) \ \ \mbox{when} \ s\leq c_2\phi(K(d)).
 \lbl{3.4}\ees

By virtue of (\ref{b.0}), it can be deduced that
 \bess
 \lim_{s\to\infty}\frac{f_2(s)}{f_1(s)}=
 \lim_{s\to\infty}\frac{\lk(k\circ K^{-1}\circ\phi^{-1}(c^{-1}_2s)\rr)^p}
 {\lk(k\circ K^{-1}\circ\phi^{-1}(c^{-1}_1s)\rr)^p}
  =\lk(\frac{c_1}{c_2}\rr)^{p\frac{1-\ell}{1-r}}.
 \eess
There is a constant $s_0>0$ such that
 \bes
 \frac 12\lk(c_1/c_2\rr)^{p\frac{1-\ell}{1-r}}f_1(s)\leq f_2(s)\leq 2\lk(c_1/c_2\rr)^{p\frac{1-\ell}{1-r}}f_1(s), \ \ \ \forall \
 s\geq s_0.\lbl{3.5}
 \ees

Let $v(t)$ be the unique positive solution of
  \bes
 v'=-f_2(v), \ \ \ t>0; \ \ \ v(0)=\infty,\lbl{3.6}
 \ees
and take $\tau=\min\{v(T),\, c_1\inf_{\Omega}\phi(K(d(x)))\}>0$.
Since $f_1(s)$ and $f_2(s)$ are positive and continuous in $(0, \infty)$, there are positive constants $C_i$ such that
$C_1f_1(s)\leq\ f_2(s)\leq C_2f_1(s)$ for all $\tau\leq s\leq s_0$.
This combined with (\ref{3.5}) yields the existence of a positive constant $C$ such that
 \bes
 C^{-1} f_1(s)\leq f_2(s)\leq Cf_1(s), \ \ \ \forall \
 s\geq \tau.\lbl{3.7}
 \ees
Since $q>0$, by Lemma \ref{l2.8}, there are a positive, continuous and increasing function $g\in RV_{q}$ and a constant $\sigma>0$ such that $\sigma g(s)\leq f_2(s)\leq g(s)$ for all $s\geq\tau$. Hence, by (\ref{3.7})
 \bes
 C^{-1}f_1(s)\leq g(s)\leq Cf_2(s), \ \ \ \forall \
 s\geq \tau.\lbl{3.8}
 \ees

Let $\Lambda>0$ be a constant and $\overline{\mathbf{u}}(x,t)=\Lambda[v(t)+z(x)]$. Then we have
 \[\overline{\mathbf{u}}_t-\Delta_p \overline{\mathbf{u}}=-\Lambda f_2(v)-\Lambda^{p-1}k^p(d(x))f(z).\]
By (\ref{3.2}), we have $c_2^{-1} z(x)\leq\phi(K(d(x)))\leq c_1^{-1} z(x)$,
which implies $k^p(d(x))f(z(x))\leq f_2(z(x))$. It follows from (\ref{3.7}) and (\ref{3.8}) that
  \bess
 \Lambda f_2(v)+\Lambda^{p-1}k^p(d)f(z)&\leq& \Lambda f_2(v)+\Lambda^{p-1}f_2(z)\\[.1mm]
   &\leq& \lk(\Lambda+\Lambda^{p-1}\rr)g(v+z)\\[.1mm]
  &\leq&\lk(\Lambda+\Lambda^{p-1}\rr)Cf_2(v+z)\\[.1mm]
  &\leq&\lk(\Lambda+\Lambda^{p-1}\rr)Cf_1(v+z).
  \eess
Since $f_1\in RV_{q}$ and $q>\max\{1, p-1\}$, by Lemma \ref{l2.5}, we can choose $\Lambda>1$ so large that
 \[(\Lambda+\Lambda^{p-1})Cf_1(v+z)\leq \alpha_\varepsilon f_1(\Lambda(v+z))=\alpha_\varepsilon f_1(\overline{\mathbf{u}}).\]
Since $\overline{\mathbf{u}}\geq z\geq c_1\phi(K(d(x)))$, by the first inequality of (\ref{3.4}), $f_1(\overline{\mathbf{u}})\leq k^p(d(x))f(\overline{\mathbf{u}})$. Hence
 \[\overline{\mathbf{u}}_t-\Delta_p \overline{\mathbf{u}}\geq-\alpha_\varepsilon k^p(d(x))f(\overline{\mathbf{u}})\geq-b(x,t)f(\overline{\mathbf{u}}), \ \ (x,t)\in
   \Omega\times(0,T-\varepsilon].\]
Thus, we obtain (\ref{3.z1}) again.

\vskip 4pt{\bf Step 2: The existence of minimal solution.}

Let $n\geq 1$ and consider the problem
   \bes\left\{\begin{array}{ll}\smallskip
 u_t-\Delta_p u=-b(x,t)f(u),\ & (x,t)\in \Omega_T,\\
 u=n, \  & (x,t)\in \Sigma_T.
  \end{array}\right.
  \lbl{3.9}\ees
Since $0$ and $n$ are the lower and upper solutions of (\ref{3.9}), it is clear that (\ref{3.9}) has a unique positive solution $u_n(x,t)$ and $u_n(x,t)$ is non-increasing in $t$. Moreover, Proposition \ref{p2.1} guarantees that $u_n(x,t)$ is strictly increasing in $n$, that is, $u_n(x,t)<u_{n+1}(x,t)$ on $\Omega_T$.

Let $\overline{\mathbf{u}}(x,t)$ be determined by Step 1. For any fixed $n$, it is clear that  $u_n(x, t)<\overline{\mathbf{u}}(x,t)$ when $(x,t)$ is near
$\Sigma_T$. Since $\overline{\mathbf{u}}(x,t)$ satisfies (\ref{3.z1}), by Proposition \ref{p2.1} we have that
$u_n(x, t)\leq \overline{\mathbf{u}}(x,t)$ in $\Omega_{T-\varepsilon}$.
It should be noticed that, for fixed small $\varepsilon
> 0$ and any compact subset $\Omega'$ of $\Omega$, $\overline{\mathbf{u}}$ is bounded on $\Omega'\times [\varepsilon,T-\varepsilon]$. As
a consequence, by standard regularity arguments, $u_n(x, t)\to
\underline{u}(x, t)$ as $n\to\infty$ uniformly on any compact subset
of $\Omega\times(0, T)$, where $\underline{u}(x,t)$ satisfies (\ref{1.1})
in the weak sense. As $u_n(x, t)$ is non-increasing in $t$, so is $\underline{u}(x,t)$. Similar to the elliptic case, it can be easily proved that $\underline{u}(x,t)=\infty$ on $\Sigma_T$; see
e.g. \cite{CW11}. Thus, $\underline{u}(x,t)$ is a solution to
(\ref{1.1})--(\ref{1.3}); in fact, it is the minimal positive
solution. Indeed, let $u(x,t)$ be any positive
solution of (\ref{1.1})--(\ref{1.3}). We can easily apply Proposition \ref{p2.1} to
conclude that $u_n(x,t)\leq u(x,t)$ in $\Omega\times (0,T]$.
Letting $n\to\infty$ we deduce $\underline{u}(x,t)\leq u(x,t)$ in
$\Omega\times(0,T)$.

\vskip 4pt {\bf Step 3: Existence of a maximal positive solution.}

We next prove the
existence of a maximal positive solution of (\ref{1.1})--(\ref{1.3}).
To achieve this, for any small $\varepsilon
> 0$, we define $\Omega_\varepsilon =\{x \in\Omega: d(x, \partial\Omega) >\varepsilon\}$. Obviously, for small
$\varepsilon, \partial\Omega_\varepsilon$ has the same smoothness as
$\partial\Omega$. We consider the following problem:
  \bes\left\{\begin{array}{ll}\smallskip
 u_t-\Delta_p u=-b(x,t)f(u),\ &
 x\in\Omega_\varepsilon\times(\varepsilon,T),\\
 u=\infty, \  & x\in
 \partial\Omega_\varepsilon\times(\varepsilon,T)\cup\bar\Omega_\varepsilon\times\{\varepsilon\}.
  \end{array}\right.
  \lbl{3.10}\ees
Let us denote by $\underline{u}^\varepsilon$ the minimal positive
solution of (\ref{3.10}). Proposition \ref{p2.1} guarantees that
$\underline{u}^{\varepsilon_1}\geq \underline{u}^{\varepsilon_2}\geq
\underline{u}$ in $\Omega_{\varepsilon_1}\times(\varepsilon_1,T)$ when $\varepsilon_1>\varepsilon_2>0$. Therefore, one
can construct a decreasing sequence $\varepsilon_n$ satisfying
$\varepsilon_n\to0$, such that $\underline{u}^{\varepsilon_n}\to\bar
u$ as $\varepsilon_n\to0$ and $\bar u$ solves
(\ref{1.1})--(\ref{1.3}). We further observe that $\bar u$ is in fact the maximal positive
solution. Indeed, for any positive
solution $u$ of (\ref{1.1})--(\ref{1.3}), it follows from the
comparison principle that
$\underline{u}^{\varepsilon_n}>u$ in $\Omega_{\varepsilon_n}\times
(\varepsilon_n,T)$ for each $n$. By taking $n\to\infty$ we obtain
$\bar u\geq u$.

\vskip 4pt {\bf Step 4: Proof of (\ref{1.11})}.

Let $\xi(t)$ and $v(t)$ be the unique positive
solution of (\ref{1.5}) and (\ref{3.6}) respectively. Since $f_2(s)\in RV_{q}$ and $q>1$, by Lemma \ref{l2.11}, there is a constant $C>0$ such that
$C^{-1} v(t)\leq\xi^*(t)\leq Cv(t)$,
here $\xi^*(t)$ is the unique positive solution of (\ref{1.z1}).

For any $0<\delta\ll1$, denote $\Omega_\delta=\{x \in\Omega: d(x, \partial\Omega) >\delta\}$. Let $z_\delta(x)$
be, respectively, the unique positive solution of
   \bessd\Delta_p z= k^p(d(x,\partial \Omega_\delta))f(z), \ \ &x\in\Omega_\delta,\\
  z=\infty, \ \ &x\in\partial\Omega_\delta\eessd
when $k$ is non-decreasing, and the unique positive solution of
  \xxx
  \Delta_p z= k^p(d(x))f(z), \ \ &x\in\Omega_\delta,\\
  z=\infty, \ \ &x\in \partial\Omega_\delta
  \lbl{3.za}\zzz
when $k$ is non-increasing (see Theorem \ref{t2.2}(i), here we emphasize that for problem (\ref{3.za}), the corresponding $k(t)=1$). Set $\xi_\delta
(t)=\xi(t-\delta)$ and $v_\delta (t)=v(t-\delta)$. From the discussion of Step 1, we can find a constant $\Lambda\geq1$, which is independent of $\delta$, such that the function
 \[\mathbf{u}^\delta(x,t)=\left\{\begin{array}{ll}\smallskip \Lambda(\xi_\delta(t)+z_\delta(x)) & {\rm if}\   k\ \mbox{\rm is\ non-increasing},\\[1mm]
   \Lambda(v_\delta(t)+z_\delta(x))\ \ &{\rm  if}\ k\ \mbox{\rm is\ non-decreasing}
   \end{array}\right.\]
satisfies
 \[\mathbf{u}^\delta _t-\Delta_p \mathbf{u}^\delta \geq-b(x,t)f(\mathbf{u}^\delta), \ \ (x,t)\in
  \Omega_\delta\times(\delta,T-\varepsilon].\]
It follows from the comparison principle that
 \bess
 \bar{u}(x,t)\leq\mathbf{u}^\delta(x,t)=\left\{\begin{array}{ll}\smallskip \Lambda(\xi_\delta(t)+z_\delta(x)) & {\rm if}\   k\ \mbox{\rm is\ non-increasing},\\[1mm]
 \Lambda(v_\delta(t)+z_\delta(x))\ \ &{\rm  if}\ k\ \mbox{\rm is\ non-decreasing}
 \end{array}\right.
 \eess
as $(x,t)\in\Omega_\delta\times(\delta,T-\varepsilon]$. Letting $\delta\to0$, and using the easily
proved fact that  $z_\delta \to z, \, \xi_\delta \to \xi$ and $v_\delta \to v$, we deduce
  \bess
  \bar{u}(x,t)\leq\left\{\begin{array}{ll}\smallskip \Lambda(\xi(t)+z(x)) &{\rm if}\
 k\ \mbox{\rm is\ non-increasing},\\[1mm]
 \Lambda(v(t)+z(x)) \ \ & {\rm if}\ k\ \mbox{\rm is\ non-decreasing}
 \end{array}\right.
 \eess
as $(x,t)\in\Omega \times(0,T-\varepsilon]$, where $z(x)$ is the unique positive solution of (\ref{3.1}) and satisfies $z(x)\leq c_2\phi( K(d(x)))$.
Thanks to $v(t)\leq C\xi^*(t)$, we conclude that (\ref{1.11}) holds.

\vskip 4pt {\bf Step 5: Proof of (\ref{1.12})}.

Case (1):  {\it $k$ is non-increasing.}

Choose $\hat\alpha_\varepsilon>0$ such that $\alpha_2(t)\leq \hat\alpha_\varepsilon$
on $[0,T-\varepsilon]$. Let $w(x)$ be the unique positive solution of
 \bessd\Delta_p w=\hat\alpha_\varepsilon k^p(d(x))f(w), \ \ &x\in\Omega,\\
  w=\infty,\ \ &x\in \partial\Omega.\eessd
Then there exist positive constants $d_1$ and $d_2$ such that
 \bes
  d_1\phi(K(d(x)))\leq w(x)\leq d_2 \phi(K(d(x))), \ \ \ x\in\Omega.
  \lbl{3.11}\ees
Following the arguments of Step 1, we have that, for any $x\in\Omega$,
 \bess\left\{\begin{array}{ll}
 k(d)\leq k\circ K^{-1}\circ\phi^{-1}(d^{-1}_1s)\ \ \ \mbox{if} \ \ \phi(K(d))\leq d_1^{-1} s,\\[2mm]
 k(d)\geq k\circ K^{-1}\circ\phi^{-1}(d^{-1}_2s)\ \ \ \mbox{if} \ \ \phi(K(d))\geq d_2^{-1} s,\end{array}\right.\eess
and the functions
 \[f^*_i(s):=\lk(k\circ K^{-1}\circ\phi^{-1}(d^{-1}_is)\rr)^pf(s)\]
satisfy $f^*_i(s)\in RV_{q}$ and
 \bes
 f^*_1(s)\geq k^p(d)f(s) \ \ \mbox{if} \ s\geq d_1\phi(K(d)), \ \ \ \ f^*_2(s)\leq k^p(d)f(s) \ \ \ \mbox{if} \ s\leq d_2\phi(K(d)).
 \lbl{3.12}\ees

Let $\eta(t)$ be the unique positive solution of
  \[\eta'=-f^*_2(\eta), \ \ \ t>0; \ \ \ \eta(0)=\infty. \]
 By Lemma \ref{l2.11}, there is a constant $C>0$ such that
 \bes
 C^{-1} \eta(t)\leq\xi^*(t)\leq C\eta(t).
 \lbl{3.13}\ees
For $0<\sigma\ll 1$, take $\tau=\min\{\eta(T+\sigma),\, d_1\inf_{\Omega}\phi(K(d(x)))\}>0$. Since $f^*_2(s)\in RV_{q}$ and $q>1$, by Lemma \ref{l2.8}, there are a positive, continuous and increasing function $g\in RV_{q}$ and a constant $\sigma>0$ such that
$\sigma g(u)\leq f^*_2(u)\leq g(u)$ for all $u\geq\tau>0$.
By Lemma \ref{l2.7}, there is a constant $c>0$, such that
$g(u)+g(v)>cg(u+v)$ for all $u,v\geq\tau$.
Similar to the discussion of Step 1, there is a constant $C>0$ such that
 \bes
 \dd C^{-1} f^*_2(u)\leq f^*_1(u)\leq Cf^*_2(u), \ \ \ \forall \
 u\geq \tau.\lbl{3.14}
 \ees

Let $\underline{\mathbf{u}}(x,t)=\kappa(\eta^\sigma(t)+w(x))$,
where $0<\kappa\ll 1$ will be chosen later and $\eta^\sigma(t)=\eta(\sigma+t)$. Since $\eta^\sigma+w\geq w\geq d_1\phi(K(d))\geq\tau$ and $g(u)\geq f^*_2(u)$ when $u\geq\tau$, by the first inequality of (\ref{3.12}) and (\ref{3.14}), there is a constant $\tilde c>0$ such that
$g(\eta^\sigma+w)\geq \tilde ck^p(d)f(\eta^\sigma+w)$ in $\Omega\times[0,T]$.
Noting that $\eta^\sigma\geq\tau$ in $[0,T]$ and $\tau\leq d_1\phi(K(d))\leq w\leq d_2\phi(K(d))$ in $\Omega$, here $d=d(x)$, we have
  \bes
  \Delta_p \underline{\mathbf{u}}-\underline{\mathbf{u}}_t&=&\kappa f^*_2(\eta^\sigma)+\kappa^{p-1}\hat\alpha_\varepsilon k^p(d)f(w)\nonumber\\
   &\geq& \kappa f^*_2(\eta^\sigma)+\kappa^{p-1}\hat\alpha_\varepsilon f^*_2(w)\nonumber\\
  &\geq&\sigma\min\{\kappa,\,\kappa^{p-1}\hat\alpha_\varepsilon\}\big[g(\eta^\sigma)+g(w)\big] \nonumber\\
  &\geq& c\sigma\min\{\kappa,\,\kappa^{p-1}\hat\alpha_\varepsilon\}g(\eta^\sigma+w)\nonumber\\
   &\geq&c\tilde c\sigma\min\{\kappa,\,\kappa^{p-1}\hat\alpha_\varepsilon\}k^p(d)f(\eta^\sigma+w)\nonumber\\
  &\geq&c\tilde c\sigma\min\{\kappa,\,\kappa^{p-1}\hat\alpha_\varepsilon\}\kappa^{-l}k^p(d)f(\underline{\mathbf{u}}) \ \ \mbox{by\ condition\ (C)}.\lbl{3.15}
  \ees
Since $l>\max\{1,p-1\}$, by (\ref{3.15}),
there is a $0<\kappa\ll 1$ such that
$c\tilde c\sigma\min\{\kappa,\,\kappa^{p-1}\hat\alpha_\varepsilon\}\kappa^{-l}\geq \hat\alpha_\varepsilon$. Consequently,
  \bes
  \underline{\mathbf{u}}_t-\Delta_p \underline{\mathbf{u}} \leq-\hat\alpha_\varepsilon k^p(d(x))f(\underline{\mathbf{u}})\leq-b(x,t)f(\underline{\mathbf{u}}), \ \ \ (x,t)\in\Omega\times(0,T-\ep].
  \lbl{3.15z}\ees

For any given $n\geq 1$, by a standard argument (see \cite{CW11}), the problem
  \bessd\Delta_p w=\hat\alpha_\varepsilon k^p(d(x))f(w), \ \ &x\in\Omega,\\
  w=n,\ \ &x\in \partial\Omega\eessd
has a unique positive solution $w_n$, and $w_n\to w$ locally uniformly in
$\Omega$ as $n\to \infty$. Since $f(s)>0$ for $s>0$, the maximum principle implies that $w_n\leq n$ on $\bar\Omega$. It follows that $w_n$ is a lower solution of (\ref{3.9}). Therefore, $u_n\geq w_n$ in $\Omega_T$ for all $n\geq 1$, and hence $\underline{u}\geq w$ in $\Omega_T$. We may assume that the constant $\kappa$, as determined above, satisfies $0<\kappa<1/2$. Then  $\underline{u}>2\kappa w$ in $\Omega_T$. Thus
 \[\underline{\mathbf{u}}-\underline{u}=\kappa(\eta^\sigma+w)-\underline{u}<
 \kappa\eta^\sigma-\frac 12\underline{u}, \ \ (x,t)\in\Omega_T.\]
Therefore,
$\dd\limsup_{(x,t)\to \Sigma_T}[\underline{\mathbf{u}}(x,t)-\underline{u}(x,t)]<0$.
Since $\underline{\mathbf{u}}$ satisfies (\ref{3.15z}), by the comparison principle,
$\underline{u}\geq \underline{\mathbf{u}}=\kappa(\eta^\sigma+w)$ in $\Omega_{T-\varepsilon}$. Taking $\sigma\to 0$ yields
$\underline{u}\geq \kappa(\eta+w)$ in $\Omega_{T-\varepsilon}$.
Thanks to (\ref{3.11}) and (\ref{3.13}), and the arbitrariness of $\varepsilon>0$, we obtain the first inequality of (\ref{1.12}).

 Case (2):  {\it $k$ is non-decreasing.}

 For any small $\sigma>0$, we consider the following auxiliary
 problems:
  \bes
 &&\xi'=-f(\xi),\ \ t>-\sigma;\ \ \ \xi(-\sigma)=\infty,
  \lbl{3.16}\\[1mm]
 &&\left\{\begin{array}{ll}
 \Delta_p z= k^p(d(x,\partial D_\sigma))f(z),\ \ &x\in D_\sigma,\\[2mm]
   z=\infty,\ &x\in\partial D_\sigma,
   \end{array}\right.\lbl{3.17}\ees
 where $D_\sigma:=\{x\in\mathbb{R}^N,d(x,\Omega)<\sigma\}$. We can
choose $\sigma$ sufficiently small such that $\partial D_\sigma$ has the same smoothness
as $\partial\Omega$. Denote by $\xi^\sigma$ and $z^\sigma$ the solutions of (\ref{3.16}) and (\ref{3.17}) respectively. It is easy to see that
for any $t\in[0, T]$ and $x\in\Omega$, both $\xi^\sigma$ and $z^\sigma$ are decreasing in $\sigma$. Hence,
 \[\tau=\min\lk\{\inf_{0<\sigma\ll 1}\xi^\sigma(T), \ \inf_{x\in\Omega,\, 0<\sigma\ll 1}z^\sigma(x)\rr\}>0.\]
As $f\in RV_{\rho}$ and $f$ is increasing, there is a constant $c>0$, such that $f(u)+f(v)>cf(u+v)$ for all $u,v\geq\tau>0$. Since $k$ is non-decreasing, we have
$\tilde ck^p(d(x))\leq1$ for some $\tilde c>0$. Set
$\underline{\mathbf{u}}(x,t)=\kappa\big[\xi^\sigma(t)+z^\sigma(x)\big]$,
where $\kappa>0$ is to be determined.
Noting that $d(x,\partial D_\sigma)>d(x,\partial\Omega)$ for all $x\in\Omega$, we have  \bess
  \underline{\mathbf{u}}_t-\Delta_p \underline{\mathbf{u}} &= &-\kappa f(\xi^\sigma)-\kappa^{p-1} k^p(d(x,\partial
  D_\delta))f(z^\sigma)\\
   &<&-\min\{\tilde c\kappa,\,\kappa^{p-1}\}k^p(d(x))\big[f(\xi^\sigma)+f(z^\sigma)\big]\\
   &\leq&-c\min\{\tilde c\kappa,\,\kappa^{p-1}\}k^p(d(x))f(\xi^\sigma +z^\sigma),\ \ \ (x,t)\in  \Omega_{T-\varepsilon}.
  \eess
Similar to Case (1), there exists a suitably small $\kappa>0$ such that
  \[-c\min\{\tilde c\kappa,\,\kappa^{p-1}\}k^p(d(x))f(\xi^\sigma +z^\sigma)\leq-\hat\alpha_\varepsilon k^p(d(x))f(\underline{\mathbf{u}}) \leq-b(x,t)f(\underline{\mathbf{u}}),\ \ \ (x,t)\in  \Omega_{T-\varepsilon}.\]
By the comparison principle
  \bes
  \underline{u}(x,t)\geq \underline{\mathbf{u}}(x,t)=\kappa\big[\xi^\sigma(t)+z^\sigma(x)\big], \ \ \ (x,t)\in \Omega_{T-\varepsilon}.
  \lbl{3.18}\ees
Clearly, $\xi^\sigma(t)\to \xi(t)$ locally uniformly on $(0,T]$ as
$\sigma\to 0^+$ and $\xi(t)$ is the unique solution of (\ref{1.5}).
Similarly, $z^\sigma(x)\to z(x)$ locally uniformly on any compact subset of $\Omega$ as $\sigma\to 0^+$, and $ z(x)$ is the unique positive solution of (\ref{3.1}). Letting $\sigma\to 0^+$ in (\ref{3.18}), and using (\ref{3.2}), the desired result is obtained since $\ep>0$ is arbitrary.
 \ \ \ \fbox{}

\section{Asymptotic behavior and uniqueness}
\setcounter{equation}{0}

In this section, we prove Theorems \ref{t1.2} and \ref{t1.3}. We first prove a lemma. Since $\rho>p-1-(p-2)/\ell$, it is easy to check that $p(1-\ell)/(r-1)<\rho-1$.

\begin{lem}\lbl{l4.1} \ For any given constant $\varsigma>0$ where $p(1-\ell)/(r-1)<\varsigma<\rho-1$, we have
 \bes
 \lim_{s\to 0^+}\frac{\phi^{-\varsigma}(K(s))}{k^p(s)}=0.
 \lbl{4.1}\ees
 \end{lem}

{\bf Proof} \ We recall that $\phi\in NRVZ_{1-r}$ (Lemma \ref{l2.2}), $K\in RVZ_{1/\ell}$ and $k\in RVZ_{(1-\ell)/\ell}$ (Lemma \ref{r2.3a}). In view of Lemma \ref{l2.4},
 \bess
 \phi^{-\varsigma}(K(s))\in RVZ_{\varsigma(r-1)/\ell}, \ \ \
 k^p(s)\in RVZ_{p(1-\ell)/\ell}, \ \ \ \frac{\phi^{-\varsigma}(K(s))}{k^p(s)}
 \in RVZ_{\varsigma(r-1)/\ell-p(1-\ell)/\ell}.
 \eess
Since $\varsigma>p(1-\ell)/(r-1)$, i.e., $\varsigma(r-1)/\ell-p(1-\ell)/\ell>0$, it is obvious that (\ref{4.1}) holds.\ \ \ \ \fbox{}

\vskip 12pt {\bf Proof of  Theorem \ref{t1.2}(i)} ~
Fix $y\in\partial\Omega$ and
$t_0 \in(0, T)$, and let $\beta_0:=\beta(y,t_0)$.  For any given small $\varepsilon\in(0,\beta_0/2)$, one can find a sufficiently small constant
$\delta\in (0,t_0)$ such that, for $(x, t)\in \Omega_T$
satisfying $|x-y|<\delta$ and $|t-t_0|<\delta$, we have
 \[\beta_0-\varepsilon\leq\frac{b(x,t)}{k^p(d(x))}\leq \beta_0+\varepsilon.\]

\vskip 4pt

{\bf Step 1}: \ We first prove the upper bound estimate
  \bes
 \limsup_{\Omega\ni x\to y}\frac{u(x,t_0)}{\phi (K(d(x)))}\leq\left(\frac{r+\ell-1}{r(\beta_0-2\varepsilon)}\right)^{\frac{r-1}{p}}.
 \lbl{4.2}
 \ees
Let $\eta(t)$ be the unique positive solution of
  \bessd
  \eta'(t)=-a f(\eta),\ \ t\in(t_0-\delta,t_0],\\
  \eta(t_0-\delta)=\infty, \ \ \eta(t_0)=1,\eessd
where $a=\frac 1\delta\int_{1}^\infty
\frac{{\rm d}s}{f(s)}>0$. Then $\eta(t)\geq 1$ on $(t_0-\delta,t_0]$. Let $\varsigma$ be given by Lemma \ref{l4.1}. Since $f$ is increasing, $f\in RV_{\rho}$ and $\rho>\varsigma+1$, by Lemma \ref{l2.6}, there is a constant $\Lambda^*>0$ such that
 \bes
 a\Lambda^{\varsigma+1}f(\eta(t))<\varepsilon f(\Lambda\eta(t)), \ \ \forall \ \Lambda\geq \Lambda^*, \ t\in (t_0-\delta,t_0].
 \lbl{4.3}\ees

Let $w_\ep(x)$ be the unique positive solution of
 \bessd
 \Delta_p w_\ep=(\beta_0-2\varepsilon)k^p(d(x))f(w_\ep), \ \ &x\in\Omega,\\
 w_\ep=\infty,\ \ &x\in\partial\Omega.\eessd
By Theorem \ref{t2.2}, there are two positive constants $d_1$ and $d_2$ such that
 \bes
 d_1\phi(K(d(x)))\leq w_\ep(x)\leq d_2\phi(K(d(x))), \ \ \ x\in\Omega.
 \lbl{4.o1}\ees
Therefore $\dd\lim_{\gamma\to 0^+}\inf_{\Omega\cap B_{2\gamma}(y)}w_\ep(x)=\infty$.
There is a constant $\gamma_0$ with $0<\gamma_0\leq\delta$ such that
  \bes
  w_\ep(x)>\Lambda^*, \ \ \ \forall \ x\in\Omega\cap B_{2\gamma}(y), \ 0<\gamma\leq\gamma_0.
  \lbl{4.4}\ees
For any fixed $0<\gamma\leq\gamma_0$, let $D\subset\Omega\cap B_{2\gamma}(y)$ be a smooth domain such
that $\partial D$ and $\partial\Omega$ coincide inside
$B_\gamma(y)$. Let $v_\ep(x)$ be a positive solution of
  \bessd
  \Delta_p v_\ep=(\beta_0-2\varepsilon)k^p(d(x))f(v_\ep), \ \ &x\in D,\\
  v_\ep=\infty,\ \ &x\in\partial D.\eessd
We note that since $d(x)=d(x,\partial\Omega)$, which may not be equal to $d(x,\partial D)$, the positive solution of the above problem may not be unique.
As $D\subset\Omega\cap B_{2\gamma}(y)$, by the comparison principle we have $v_\ep(x)\geq w_\ep(x)$ in $D$. Hence, by (\ref{4.4}),
  \bes
  v_\ep(x)\geq w_\ep(x)>\Lambda^*, \ \ \ \forall \ x\in D.
  \lbl{4.6}\ees
From the choice of $D$, it is clear that $d(x)=d(x,\partial D)$ when $x\in D$ and is near $y$. Evidently,
 \[\lim_{D\ni x\to y}\frac{(\beta_0-2\varepsilon)k^p(d(x))}
  {k^p(d(x,\partial D))}=\beta_0-2\varepsilon.\]
In view of Remark \ref{r2.2} we have
  \bes
 \lim_{D\ni x\to y}\frac{v_\ep(x)}{\phi(K(d(x)))}
 =\lim_{D\ni x\to y}\frac{v_\ep(x)}{\phi(K(d(x,\partial D)))}
 =\left(\frac{r+\ell-1}{r(\beta_0-2\varepsilon)}\right)^{\frac{r-1}{p}}.
 \lbl{4.5}\ees

We now consider $\Omega_\sigma:=\{x\in\Omega: d(x)\geq\sigma\}$ for
sufficiently small $\sigma\in[0,\gamma/2)$. For each such
$\Omega_\sigma$, we can construct a
smooth domain $D_\sigma\subset\Omega_\sigma\cap
B_{2\gamma}(y)\subset D$ such that $\partial D_\sigma$ and
$\partial\Omega_\sigma$ coincide inside $B_{\gamma}(y)$, and
$D_\sigma$ varies continuously with $\sigma$ for all small
non-negative $\sigma$. We may also require that $D_\sigma\subset
D_{\sigma'}$ when $\sigma>\sigma'$ and $D_\sigma\to D$ as $\sigma\to 0^+$. By Theorem \ref{t2.2}, the problem
  \bessd\Delta_p v_\sigma=(\beta_0-2\varepsilon)k^p(d(x))f(v_\sigma), \ \ &x\in D_\sigma,\\ v_\sigma=\infty,\ \ &x\in\partial D_\sigma\eessd
has a unique positive solution, denoted by $v_\sigma$ (with $k(t)=1$, $\beta(y)=(\beta_0-2\varepsilon)k^p(d(y))$ for $y\in\partial D_\sigma$).
Applying the comparison principle and (\ref{4.6}) we get
$v_\sigma(x)\geq v_\ep(x)\geq w_\ep(x)>\Lambda^*$ in $D_\sigma$. By further using the elliptic regularity, we see that $v_\sigma$ decreases to $v_\ep$ as $\sigma$ decreases to $0$.

Set $u_\sigma(x,t)=\eta(t)v_\sigma(x)$. Then for $(x,t)\in D_\sigma\times (t_0-\delta,t_0]$,
  \bes
 (u_\sigma)_t-\Delta_p u_\sigma &=&\eta'v_\sigma-\eta^{p-1}\Delta_p v_\sigma\nonumber\\
  &=&-a v_\sigma f(\eta)-(\beta_0-2\varepsilon)k^p(d(x))\eta^{p-1}(t)f(v_\sigma)\nonumber\\
   &=&-a v^{-\varsigma}_\sigma v^{\varsigma+1}_\sigma f(\eta)-
  (\beta_0-2\varepsilon)k^p(d(x))\eta^{p-1}f(v_\sigma).
  \lbl{4.7}
  \ees
Thanks to the facts that $f(s)/s^{p-1}$ is increasing and $\eta\geq 1$, one has
  \bes
  \eta^{p-1}f(v_\sigma)\leq f(\eta v_\sigma)=f(u_\sigma)),\ \ (x,t)\in D_\sigma\times (t_0-\delta,t_0].
  \lbl{4.8}\ees
Noting that $v_\sigma(x)>\Lambda^*$ and $\varsigma>0$, and taking into account (\ref{4.3}), we have
 \bes
 av_\sigma^{\varsigma+1}f(\eta)<\varepsilon f(v_\sigma\eta)
   =\varepsilon f(u_\sigma),\ \ (x,t)\in D_\sigma\times (t_0-\delta,t_0].
 \lbl{4.9}\ees
If $k(0)>0$, then $k^p(d(x))$ has a positive lower bounded in $\Omega$. In view of $\dd\lim_{\gamma\to 0^+}\dd\inf_{\Omega\cap B_{2\gamma}(y)}w_\ep=\infty$, we can choose $\gamma$ small enough such that $w_\ep^{-\varsigma}(x)<k^p(d(x))$. Hence
 \bes
 v^{-\varsigma}_\sigma(x)<k^p(d(x)), \ \ \ \forall \ x\in D_\sigma.
 \lbl{4.10}\ees
If $k(0)=0$, by Lemma \ref{l4.1},
 \bes
 \lim_{d(x)\to 0^+}\frac{\phi^{-\varsigma}(K(d(x)))}{k^p(d(x))}=0.
 \lbl{4.11}\ees
In view of $v_\sigma(x)\geq w_\ep(x)$ in $D_\sigma$ and the estimates (\ref{4.o1}), it follows that $v_\sigma(x)\geq w_\ep(x)\geq d_1\phi(K(d(x)))$ in $D_\sigma$. Therefore
 \bes
 \frac{v_\sigma^{-\varsigma}(x)}{k^p(d(x))}\leq
 d_1^{-\varsigma}\frac{\phi^{-\varsigma}(K(d(x)))}{k^p(d(x))},\ \ \ \forall \ x\in D_\sigma.\lbl{4.12}\ees
It is clear that $d(x)\to 0^+$ holds uniformly on $\overline{D}_\sigma$ as $\gamma\to 0^+$. By virtue of (\ref{4.11}) and (\ref{4.12}), one can choose $\gamma$ small enough such that (\ref{4.10}) is true.

It follows from (\ref{4.7})--(\ref{4.11}) that ($d=d(x)$)
  \bes
  (u_\sigma)_t-\Delta_p u_\sigma
     &=&-av^{-\varsigma}_\sigma v^{\varsigma+1}_\sigma f(\eta)-
  (\beta_0-2\varepsilon)k^p(d)\eta^{p-1}f(v_\sigma)\nonumber\\
  &\geq&-(\beta_0-\varepsilon)k^p(d)f(u_\sigma)\nonumber\\
  &\geq&-bf(u_\sigma),\ \ \ (x,t)\in D_\sigma\times (t_0-\delta,t_0].
  \lbl{4.13}
  \ees
It is obvious that
 \bess
 &u(x,t_0-\delta)<u_\sigma(x,t_0-\delta), \ \ \ x\in D_\sigma, &\\ &u(x,t)|_{\partial D_\sigma}<u_\sigma(x,t)|_{\partial D_\sigma}, \ \ \ t\in(t_0-\delta,t_0].&
 \eess
By (\ref{4.13}) and the comparison principle,
$u(x,t)\leq u_\sigma(x,t)=\eta(t)v_\sigma(x)$ in $D_\sigma\times(t_0-\delta,t_0]$.
Letting $\sigma\to0$, one has
$u(x,t)\leq \eta (t)v_\ep(x)$ in $D\times(t_0-\delta,t_0]$. Hence, by (\ref{4.5}),
 \bess
 \limsup_{\Omega\ni x\to y}\frac{u(x,t_0)}{\phi (K(d(x)))}\leq \limsup_{D\ni x\to y}\frac{v_\ep(x)}{\phi( K(d(x)))}
 =\left(\frac{r+\ell-1}{r(\beta_0-2\varepsilon)}\right)^{\frac{r-1}{p}}.
 \eess
We thus obtain the estimate (\ref{4.2}).

\vskip 4pt

{\bf Step 2}: \ Now we establish the lower bound estimate
  \bes
 \liminf_{\Omega\ni x\to y}\frac{u(x,t_0)}{\phi (K(d(x)))}\geq\left(\frac{r+\ell-1}{r(\beta_0+\varepsilon)}\right)^{\frac{r-1}{p}}.
 \lbl{4.14}
 \ees

Choose a constant $A_0>0$ such that $\alpha_2(t)\leq A_0$ in $[0,t_0]$. Let $w$ and $z$ be the unique positive solutions of following problems, respectively:
 \bess
 &\Delta_p w=A_0k^p(d(x))f(w), \ \ x\in\Omega;\ \ \ w=\infty,\ \ x\in\partial\Omega,&\\
 &\Delta_p z=(\beta_0+\varepsilon)k^p(d(x))f(z), \ \ x\in\Omega;\ \
 \ z=\infty,\ \ x\in\partial\Omega.&
  \eess
According to Theorem \ref{t2.2}, there is a constant $C>0$ such that
 \bes
 C^{-1}w(x)\leq z(x)\leq Cw(x), \ \ \ \forall \ x\in\Omega.\lbl{4.n1}\ees
Let $w_n$ be the unique positive solution of
 \bessd\Delta_p w=A_0k^p(d(x))f(w), \ \ &x\in\Omega,\\
 w=n,\ \ &x\in\partial\Omega.\eessd
Then $w_n$ is increasing in $n$ and $w_n\to w$ uniformly on any compact subset of $\Omega$.
Thanks to $b(x,t)\leq\alpha_2(t)k^p(d(x))\leq A_0k^p(d(x))$ for all $(x,t)\in\Omega\times[0,t_0]$, we see that $w_n$ satisfies
 \bessd\Delta_p w_n=A_0k^p(d(x))f(w_n)\geq b(x,t)f(w_n), \ \ &x\in\Omega,\\
  w_n=n,\ \ &x\in\partial\Omega.\eessd
Proposition \ref{p2.1} asserts that $w_n\leq u$ in $\Omega\times[0,t_0]$ for all $n$. Hence, $w\leq u$ in $\Omega\times[0,t_0]$. This combines with (\ref{4.n1}) to yield
$z\leq Cw\leq Cu$ in $\Omega\times[0,t_0]$. Denote $\zeta_0=C^{-1}$. Then we have
 \bes
 \zeta_0 z(x)\leq u(x,t_0-\delta), \ \ x\in\Omega.\lbl{4.15}
  \ees
Let $\zeta$ be the unique positive solution of
  \bessd\zeta'(t)=\tilde a f(\zeta), \ \ t\in(t_0-\delta,t_0],\\
  \zeta(t_0-\delta)=\zeta_0, \ \ \zeta(t_0)=1,\eessd
where $\tilde a=\frac 1\delta\int_{\zeta_0}^1
\frac{{\rm d}s}{f(s)}>0$. Then $\zeta(t)\leq 1$ on $[t_0-\delta,t_0]$.

\vskip 4pt

We first consider the case $1<p<2$.
As above, let $D\subset\Omega\cap B_{2\delta}(y)$ be a smooth domain such
that $\partial D$ and $\partial\Omega$ coincide inside
$B_\delta(y)$. Take $\psi=\frac 12z|_{\partial D}$,
and let $\lk\{\psi_n\rr\}_{n=1}^\infty$ be an increasing sequence of non-negative smooth functions defined on $\partial D$ with the property that
 \begin{quote}
  \qquad $\psi_n\big|_{\overline{\partial D\cap\partial\Omega}}=n$ \ and \ $\psi_n\to\psi$ uniformly on any compact subset of $\partial D\setminus\overline{\partial D\cap\partial\Omega}$.
 \end{quote}
Let $A>\tilde af(1)$ be a given constant. Then for any $m\geq 1$, the problem
  \[\Delta_pv=Av+(\beta_0+\varepsilon)\min\{m,k^p(d(x)\}f(v),\ \ x\in D;\ \ \
   v|_{\partial D}=\psi_n\]
has a the unique positive solution, denoted by $v_n^m$; cf. \cite{CW11}. By the comparison principle, $v_n^m\geq v_n^{m+1}$.  Thus $v_n=
\dd\lim_{m\to\infty} v^m$ exists, and one easily sees by standard
elliptic regularity that $v_n$ is a solution to
 \bessd\Delta_pv=Av+(\beta_0+\varepsilon)k^p(d(x)f(v),\ \ &x\in D,\\
   v=\psi_n, &x\in\partial D.\eessd
The comparison principle infers that $v_n$ is unique,
$v_n \leq v_{n+1}$ in
$D$ since $\psi_n \leq \psi_{n+1}$ on $\partial D$, and $v_n \leq v^*$, where $v^*$ is the unique positive solution of
 \xxx\Delta_pv=Av+(\beta_0+\varepsilon)k^p(d(x))f(v),\ \ &x\in D,\\
   v=\infty,&x\in\partial D.
 \lbl{4.15y}\zzz
Since $p<2$ and $A>0$, the function $As/s^{p-1}=As^{2-p}$ is increasing in $s>0$, and hence the comparison principle holds for the problem (\ref{4.15y}). The existence and uniqueness of $v^*$ can be proved by the similar methods of \cite{CW11, LPW12}.

Thus $v:=\dd\lim_{n\to\infty} v_n$ exists,
and by the elliptic regularity we find that $v$ is a positive solution
of
  \xxx
  \Delta_pv=Av+(\beta_0+\varepsilon)k^p(d(x)) f(v),\ \ &x\in D,\\
   v=\psi=\frac 12 z,&x\in\partial D.
  \lbl{4.15z}\zzz
In fact, by the interior regularity it is easy to show that  $v$ satisfies the differential equation of (\ref{4.15z}). Using the boundary estimate we can prove that $v$ is continuous in $\overline D\setminus\overline{\partial D\cap\partial\Omega}$. Hence, $v=\psi$ in $\partial D\setminus\overline{\partial D\cap\partial\Omega}$ in the classical sense. Now we prove that for any $x_0\in \overline{\partial D\cap\partial\Omega}$, the limit $\dd\lim_{D\ni x\to x_0}v(x)=\infty$ holds. If this is not true, then there exist $x_0\in\overline{\partial D\cap\partial\Omega}$, a sequence $\{x_l\}_{l=1}^\infty\subset D$ and a constant $M>0$, such that $x_l\to x_0$ and $v(x_l)\leq M$. Since $v_n\leq v$ for all $n$, we have $v_n(x_l)\leq M$ for all $n$ and $l$. Letting $l\to\infty$, we see that $v_n(x_0)\leq M$ for all $n$. This is a contradiction since $v_n(x_0)=\psi_n(x_0)=n$ for all $n$.

The comparison principle asserts $v\leq z$ in $D$. Because the comparison principle holds for the problem (\ref{4.15z}), similar to the proof of \cite[Theorem 1.2(i)]{CW11} or \cite[Theorem 1.2]{LPW12}, we also have that
 \bes
 \liminf_{D\ni x\to y }\frac{v(x)}{\phi (K(d(x)))}
 =\left(\frac{r+\ell-1}{r(\beta_0+\varepsilon)}\right)^{\frac{r-1}{p}}.
 \lbl{4.16}\ees

Set $u^*(x,t)=\zeta(t)v(x)$
for $(x,t)\in D\times [t_0-\delta,t_0]$. Clearly, by (\ref{4.15}),
$u^*(x,t_0-\delta)=\zeta_0v(x)\leq \zeta_0z(x)\leq u(x,t_0-\delta)$ in $D$.
It is also evident that $u^*\leq u$ on $\partial D\times
[t_0-\delta,t_0]$. A direct computation yields
 \bess
 u^*_t-\Delta_pu^*&=&\zeta'v-\zeta^{p-1}\Delta_pv\\
 &=&\tilde a f(\zeta)v-\zeta^{p-1}[Av+(\beta_0+\varepsilon) k^p(d(x))f(v)]\\[1mm]
 &=&\left(\tilde a f(\zeta)-A\zeta^{p-1}\rr)v-(\beta_0+\varepsilon) k^p(d(x))\zeta^{p-1}f(v).\eess
Thanks to the facts that $f(s)/s^{p-1}$ is increasing in $s>0$ and $\zeta(t)\leq 1$ in $[t_0-\delta, t_0]$, one has
$f(\zeta)\leq f(1)\zeta^{p-1}$ and $\zeta^{p-1}f(v)\geq f(\zeta v)=f(u^*)$.
As $A>\tilde af(1)$, it follows that
 \bess
 u^*_t-\Delta_pu^*\leq-(\beta_0+\varepsilon)k^p(d(x))f(u^*), \ \ (x,t)\in D\times [t_0-\delta,t_0].\eess
We can apply the comparison principle to conclude that $u^*\leq u$ in $D\times [t_0-\delta,t_0]$. In particular, $v(x)=u^*(x,t_0)\leq u(x,t_0)$ in $D$.
By (\ref{4.16}), it follows that
 \bess
 \liminf_{\Omega\ni x\to y }\frac{u(x,t_0)}{\phi (K(d(x)))}\geq
 \liminf_{D\ni x\to y }\frac{v(x)}{\phi (K(d(x)))}
 =\left(\frac{r+\ell-1}{r(\beta_0+\varepsilon)}\right)^{\frac{r-1}{p}}.
 \eess
Hence (\ref{4.14}) holds. The desired result clearly follows from
(\ref{4.2}) and (\ref{4.14}), since $\varepsilon > 0$ can be
arbitrarily small.

\vskip 2pt Next, we consider the case $p\geq 2$. As above, let $A>\tilde af(1)$ be a given constant. By arguments similar to those of \cite[Theorem 4.4]{DG03} and \cite[Theorem 1.2]{CW11}, it can be proved that the problem
 \bessd
 \Delta_p z=Az^{p-1}+(\beta_0+\varepsilon)k^p(d(x))f(z), \ \ &x\in\Omega,\\
 z=\infty,\ \ &x\in\partial\Omega\eessd
has a unique positive solution, denoted by $\hat z$. The comparison principle yields $\hat z\leq z$ in $\Omega$, and hence
 \bes
 \zeta_0\hat z(x)\leq\zeta_0 z(x)\leq u(x,t_0-\delta), \ \ x\in\Omega\lbl{4.19a}
  \ees
by (\ref{4.15}). Moreover, there are two positive constants $d_1'$ and $d_2'$ such that
 \[d_1'\phi(K(d(x)))\leq \hat z(x)\leq d_2'\phi(K(d(x))), \ \ \ x\in\Omega.\]
Therefore $\dd\lim_{\gamma\to 0^+}\inf_{\Omega\cap B_{2\gamma}(y)}\hat z(x)=\infty$.
There is a constant $\gamma>0$ such that
  \bes
  \hat z(x)>2, \ \ \ \forall \ x\in\Omega\cap B_{2\gamma}(y).
  \lbl{4.4zz}\ees
As above, let $\hat D\subset\Omega\cap B_{2\gamma}(y)$ be a smooth domain such
that $\partial \hat D$ and $\partial\Omega$ coincide inside
$B_\gamma(y)$. Similar to the above, the problem
  \xxx
 \Delta_pv=Av^{p-1}+(\beta_0+\varepsilon)k^p(d(x)) f(v),\ \ &x\in \hat D,\\[1mm]
 v=\frac 12\hat z, &x\in\partial \hat D
  \lbl{4.15zz}\zzz
has a positive solution, denoted by $\hat v$, and $\hat v$ satisfies (\ref{4.16}). Moreover, by the comparison principle, $\hat v\leq\hat z$ in $\hat D$. Since the function $w=\frac 12\hat z$ satisfies
 \[\Delta_pw\geq Aw^{p-1}+(\beta_0+\varepsilon)k^p(d(x)) f(w),\]
the comparison principle gives $\hat v\geq w=\frac 12\hat z$ in $\hat D$. Hence $\hat v>1$ in $\hat D$ by (\ref{4.4zz}).

Set $\hat u(x,t)=\zeta(t)\hat v(x)$
for $(x,t)\in \hat D\times [t_0-\delta,t_0]$. Then, as $\hat v>1$ in $\hat D$ and $p\geq 2$, similar to the above, we have that
 \bess
 \hat u_t-\Delta_p\hat u&=&\tilde a f(\zeta)\hat v-\zeta^{p-1}[A\hat v^{p-1}+(\beta_0+\varepsilon) k^p(d(x))f(\hat v)]\\[1mm]
 &=&\left(\tilde a f(\zeta)-A\zeta^{p-1}\hat v^{p-2}\rr)\hat v-(\beta_0+\varepsilon) k^p(d(x))\zeta^{p-1}f(\hat v)\\[1mm]
  &\leq&-(\beta_0+\varepsilon)k^p(d(x))f(\hat u), \ \ (x,t)\in \hat D\times [t_0-\delta,t_0].\eess
Since $\hat v(x)$ satisfies (\ref{4.16}), the rest of the proof is the same as that of the case $1<p<2$.  \ \ \ \fbox{}

\vskip 12pt {\bf Proof of  Theorem \ref{t1.2}(ii)} ~Let $x_0
\in\Omega$ be fixed. Then, for any given small $\varepsilon> 0$, we can find
a small ball $B_r(x_0)$ and small $t_0
> 0$ such that $\bar B_r(x_0)\subset\Omega$ and
 \[0<b_0-\varepsilon\leq b(x,t)\leq b_0+\varepsilon\]
for all $x\in B_r(x_0), t\in [0,t_0]$,  where $b_0=b(x_0,0)$.

\vskip 2pt

{\bf Step 1} \ Let $\mu^*_\varepsilon(t)$ be the unique positive solution of
  \[\big(\mu^*_\varepsilon\big)'=-(b_0-\varepsilon)f(\mu^*_\varepsilon),\ \ t>0;\ \ \
  \mu^*_\varepsilon(0)=\infty.\]
We shall prove that
 \bes
 \limsup_{t\to 0^+}\frac{u(x_0,t)}{\mu^*_\varepsilon(t)}\leq 1.
 \lbl{4.17}\ees

By Proposition \ref{p2.4}, $u^{-\rho}f(u)=M(u)\exp\left\{\int_b^u\frac{\varphi(t)}t{\rm d}t\right\}$,
where $M(u)$ satisfies $\dd\lim_{u\to\infty}M(u)=M^*>0$.
It follows that there is a sequence $\{a_n\}$ with $\dd\lim_{n\to\infty}a_n=\infty$, such that $\frac{n-1}nM^*\leq M(u)\leq\frac n{n-1}M^*$ when $u\geq a_n$. Define
 \[f_0(u)=M^*u^{\rho}\exp\left\{\int_b^u\frac{\varphi(t)}t{\rm d}t\right\}.\]
Then $f_0(u)$ satisfies
 \bes
 \frac{n-1}nf_0(u)\leq f(u)\leq\frac n{n-1}f_0(u), \ \ \ \forall \ u\geq a_n.
 \lbl{4.18}\ees
Let $\nu>0$ be small such that $1+\nu<\rho$. Note that $\rho>p-1$ and $\rho>1+\nu$, similar to the proof of Lemma \ref{l2.8}, we can prove that $f_0(u)/u^{1+\nu}$ is increasing when $u\geq A\gg 1$. Certainly, $f_0(u)/u$ is increasing for $u\geq A$, and hence
 \bes
 f_0(u)+f_0(v)\leq f_0(u+v), \ \ \ \forall \ u,v\geq A.
 \lbl{4.19}\ees
It can be assumed that $a_n\geq A$ for all $n$ without loss of generality.

Let $\eta_n(t)$, $\zeta_n(t)$, $z_n(x)$ and $\hat z(x)$ be solutions of the following problems, respectively:
  \bess
  &\eta_n'=-(b_0-\varepsilon)\dd\frac{n-1}nf_0(\eta_n),\ \ t>0;\ \ \ \eta_n(0)=\infty;&\\[1mm]
  &\zeta_n'=-(b_0-\varepsilon)\dd\frac n{n-1}f_0(\zeta_n),\ \ t>0;\ \ \ \zeta_n(0)=\infty;&\\[1mm]
  &\Delta_p z_n=(b_0-\varepsilon)\dd\frac{n-1}nf_0(z_n),\ \ x\in B_r(x_0);\ \ \ z_n=\infty, \ \ x\in \partial B_r(x_0);&\\[1mm]
  &\Delta_p\hat z=(b_0-\varepsilon)f_0(\hat z), \ \ x\in B_r(x_0);\ \ \
  \hat z=\infty \ \ x\in\partial B_r(x_0).&
  \eess
By a simple comparison argument, we have $z_n(x)\geq\hat z(x)$ in $B_r(x_0)$. Similar to the proof of Theorem 6.1 in \cite{Du}, we can prove that $\dd\lim_{r\to 0^+}\min_{B_r(x_0)}\hat z(x)=\infty$. Choosing $r$ small enough, we may assume that $\hat z(x)>A$, and hence $z_n(x)>A$ on $\overline{B_r(x_0)}$. It is also evident that there is $t_n>0$ such that $\eta_n(t)\geq a_n,\, \zeta_n(t)\geq a_n$ for all $0<t\leq t_n$.

Define $u_n(x,t)=\eta_n(t)+z_n(x)$. By (\ref{4.19}) we have that, for $(x,t)\in B_r(x_0)\times(0,t_n]$,
 \bess
 (u_n)_t-\Delta_pu_n&=&-(b_0-\varepsilon)\frac{n-1}n\big[f_0(\eta_n)+f_0(z_n)\big]\\[0.5mm]
 &\geq&-(b_0-\varepsilon)\frac{n-1}nf_0(\eta_n+z_n)\\[0.5mm]
 &\geq&-(b_0-\varepsilon)f(\eta_n+z_n)\\[0.5mm]
 &=&-(b_0-\varepsilon)f(u_n).
 \eess
By the comparison principle we obtain
  \bes
  u(x,t)\leq u_n(x,t)=\eta_n(t)+z_n(x), \ \ (x,t)\in B_r(x_0)\times(0,t_n].
  \lbl{4.20}\ees

Thanks to (\ref{4.18}), by a simple comparison argument, we have
$\zeta_n(t)\leq\mu^*_\varepsilon(t)\leq\eta_n(t)$ in $(0,t_n]$.
Let $l_n=[n/(n-1)]^{2/\nu}>1$. If we can prove $\eta_n(t)\leq l_n\zeta_n(t)$ in $(0,t_n]$, then
   \bes
 \eta_n(t)\leq l_n\zeta_n(t)\leq l_n\mu^*_\varepsilon(t), \ \ \ t\in(0,t_n].\lbl{4.21}\ees
Since $f_0(u)/u^{1+\nu}$ is increasing for $u\geq A$ and $\zeta_n(t)\geq a_n\geq A$ for $0<t\leq t_n$, it follows that $l_n^{1+\nu}f_0(\zeta_n(t))\leq f_0(l_n\zeta_n(t))$ for all $0<t\leq t_n$. Therefore
 \[(l_n\zeta_n)'=-l_n(b_0-\varepsilon)\frac n{n-1}f_0(\zeta_n)\geq
-(b_0-\varepsilon)\frac {n-1}nf_0(l_n\zeta_n)=\eta'_n, \ \ \ 0<t\leq t_n.\]
Consequently, $\eta_n(t)\leq l_n\zeta_n(t)$ in $(0,t_n]$ by the comparison principle. By (\ref{4.20}) and (\ref{4.21}), $u(x,t)\leq l_n\mu^*_\varepsilon(t)+z_n(x)$ in $B_r(x_0)\times(0,t_n]$. Hence, $\limsup\limits_{t\to 0^+}\frac{u(x_0,t)}{\mu^*_\varepsilon(t)} \leq l_n$ for all $n$. Taking $n\to\infty$, we obtain (\ref{4.17}).

Choose $0<\varepsilon_n\to 0^+$ and set $\varepsilon=\varepsilon_n$. Taking into account that $\dd\lim_{t\to 0^+}\mu^*_{\varepsilon_n}(t)=\infty$, similar to the above arguments we can prove that there are $\delta_n\to 0^+$ and $t_n'\to 0^+$ such that
  \bes
  \mu^*_{\varepsilon_n}(t)\leq(1+\delta_n)\tau(t), \ \ \forall \ t\in(0,t_n'), \ n\gg 1,
  \lbl{4.25a}\ees
where $\tau(t)$ is the unique positive solution of (\ref{1.13}). It is deduced from (\ref{4.17}) and (\ref{4.25a}) that
$\dd\limsup_{t\to 0^+}\frac{u(x_0,t)}{\tau(t)}\leq 1+\delta_n$ for all $n$.
The limit (\ref{1.15}) is obtained by take $n\to\infty$.

\vskip 4pt

{\bf Step 2} \ Let $\widetilde \mu_\varepsilon(t)$ be the positive solution of
  \[\widetilde \mu_\varepsilon'=-(b_0+\varepsilon)f(\widetilde \mu_\varepsilon),\ \ t>0;\ \ \ \widetilde  \mu_\varepsilon(0)=\infty.\]
Under the conditions that $p>2N/(N+2)$ and $f(s)/s$ is increasing for $s>0$, we shall prove that
 \bes
 \liminf_{t\to 0^+}\frac{u(x_0,t)}{\widetilde \mu_\varepsilon(t)}\geq 1.
 \lbl{4.22}\ees

We first consider the case $p\geq 2$. Let $\lambda_1$ be the first eigenvalue of the problem
  \bessd-\Delta_p\vp=\lambda_1\vp^{p-1}, \ \ &x\in B_r(x_0),\\
  \vp= 0, \ \ &x\in \partial B_r(x_0),\eessd
and $\vp$ with ${\rm sup}_{B_r(x_0)}\vp(x)= 1$ be the positive eigenfunction corresponding to $\lambda_1$. Then $\lambda_1>0$. Obviously,
$0 <\vp(x)<1$ in $B_r(x_0)\setminus\{x_0\}$ and $\vp(x_0)=1$.
Let $\mu_\varepsilon$ be the unique positive solution of
 \[\mu_\varepsilon'=-\lambda_1\mu_\varepsilon^{p-1}-(b_0+\varepsilon)f(\mu_\varepsilon), \ \ t > 0;\ \ \ \mu_\varepsilon(0) = \infty.\]
For any $\sigma:0<\sigma\ll 1$, set $\omega(x,t)=\mu_\varepsilon(\sigma+t)\vp(x)$. Since $0\leq\vp(x)\leq 1$ and $f(s)/s$ is increasing in $s>0$, it follows that
 \bess
 \omega_t-\Delta_p\omega&=&\vp\mu_\varepsilon'(\sigma+t)-
 \mu_\varepsilon^{p-1}(\sigma+t)\Delta_p\vp\\
 &=&\lambda_1\mu_\varepsilon^{p-1}(\sigma+t)(\vp^{p-1}-\vp)-(b_0+\varepsilon)\vp f(\mu_\varepsilon)\\
 &\leq&-(b_0+\varepsilon)f(\mu_\varepsilon\vp)=-(b_0+\varepsilon)f(\omega), \ \ (x,t)\in B_r(x_0)\times[0,t_0].
 \eess
Noting that $\vp=0$ on $\partial B_r(x_0)$, and $u=\infty$ on $\overline{B_r(x_0)}\times\{0\}$, the above inequality shows that $\omega(x,t)$ is a lower solution of the problem
 \bes\left\{\begin{array}{ll}\smallskip
 v_t-\Delta_pv=-b(x,t)f(v), \ \ &(x,t)\in B_r(x_0)\times(0,t_0),\\
 v=u, \ \ &(x,t)\in\partial B_r(x_0)\times(0,t_0)\cup \overline{B_r(x_0)}\times\{0\}.
 \end{array}\right.
 \lbl{4.23}\ees
Clearly $u$ solves (\ref{4.23}).
The comparison argument then implies $\mu_\varepsilon(\sigma+t)\vp(x)\leq u(x,t)$ in $B_r(x_0)\times(0,t_0]$. Letting $\sigma\to 0^+$, we deduce $\mu_\varepsilon(t)\vp(x)\leq u(x,t)$ in $B_r(x_0)\times(0,t_0]$. In
particular,
\bes
 \mu_\varepsilon(t)\leq u(x_0,t), \ \ \ t\in (0,t_0].
 \lbl{4.24}\ees

Noting that $f\in RV_{\rho}$ and $\rho>p-1$, in view of $\mu_\varepsilon(t)\to\infty$ as $t\to 0^+$, there exists $t_n$ with $0<t_n\leq t_0$ such that
 \[\frac{\lambda_1\mu_\varepsilon^{p-1}(t)}{(b_0+\varepsilon)f(\mu_\varepsilon(t))}\leq \frac 1n, \ \ \ 0<t\leq t_n.\]
Hence $\mu_\varepsilon(t)$ satisfies
\[\mu_\varepsilon'\geq-(b_0+\varepsilon)(1+1/n)f(\mu_\varepsilon), \ \ 0<t\leq t_n;\ \ \ \mu_\varepsilon(0) = \infty.\]
Using the arguments of Step 1 we can prove that there exist $\ell_n\nearrow 1$ and $t_n^*$ with $0<t_n^*\leq t_n$, such that
$\ell_n\widetilde \mu_\varepsilon(t)\leq \mu_\varepsilon(t)$, \ $\forall \ 0<t<t_n^*, \ n\gg 1$.
This combined with (\ref{4.24}) gives
$\ell_n\widetilde \mu_\varepsilon(t)\leq u(x_0,t)$, \ $\forall \ 0<t<t_n^*, \ n\gg 1$. Therefore, $\liminf\limits_{t\to 0^+}\frac{u(x_0,t)}{\widetilde \mu_\varepsilon(t)}\geq \ell_n$, \ $\forall \ n\gg 1$. Letting $n\to\infty$ we get (\ref{4.22}).

\vskip 4pt

Now we consider the case $2N/(N+2)<p<2$. By \cite[Theorems I and II]{Ota}, there is a constant $\lambda>0$ such that the problem
 \bessd
 -\Delta_p\vp=\lambda\vp, \ \ &x\in B_r(x_0),\\
  \vp= 0, \ \ &x\in \partial B_r(x_0)\eessd
has a positive solution $\vp\in W^{1,p}_0(\Omega)\cap L^\infty(\Omega)$ with $0 <\vp(x)<1$ in $B_r(x_0)\setminus\{x_0\}$ and $\vp(x_0)=1$. Let $\eta_\varepsilon$ be the unique positive solution of
 \[\eta_\varepsilon'=-\lambda\eta_\varepsilon^{p-1}-(b_0+\varepsilon)f(\eta_\varepsilon), \ \ t > 0;\ \ \ \eta_\varepsilon(0) = \infty.\]
Similar to the above we can prove that $\eta_\varepsilon(\sigma+t)\vp(x)\leq u(x,t)$ in $B_r(x_0)\times(0,t_0]$ provided that $0<\sigma\ll 1$. Letting $\sigma\to 0^+$ we deduce $\eta_\varepsilon(t)\vp(x)\leq u(x,t)$ in $B_r(x_0)\times(0,t_0]$. In
particular, $\eta_\varepsilon(t)\leq u(x_0,t)$ in $(0,t_0]$. By the same argument as above we obtain (\ref{4.22}).

As above, choose $0<\varepsilon_n\to 0^+$ and set $\varepsilon=\varepsilon_n$. Taking into account
that $\dd\lim_{t\to 0^+}\widetilde \mu_{\varepsilon_n}(t)=\infty$, similar to the above arguments we can prove that there are $\sigma_n\to 0^+$ and $t_n'\to 0^+$ such that
  \bes
  (1-\sigma_n)\tau(t)\leq\widetilde \mu_{\varepsilon_n}(t), \ \ \forall \ t\in(0,t_n'), \ n\gg 1,
  \lbl{4.25}\ees
here $\tau(t)$ is the unique positive solution of (\ref{1.13}). By (\ref{4.22}) and (\ref{4.25}), $\liminf\limits_{t\to 0^+}\frac{u(x_0,t)}{\tau(t)}\geq 1-\sigma_n$.
The limit (\ref{1.16}) is obtained by taking $n\to\infty$. The proof is complete. \ \ \ \ \fbox{}

\vskip 4pt {\bf Proof of  Theorem \ref{t1.3}} ~ Under our assumptions, it is easy to see that $f^*(u)=f(u)$, and hence $\xi^*(t)=\xi(t)$. By (\ref{1.11}) and (\ref{1.12}), there is a constant $l>1$ such that $\underline u(x,t)\leq\bar u(x,t)\leq l\underline u(x,t)$ in $\Omega_T$. The remainder of the proof is similar to that of \cite[Theorem 1.4]{DPP11}. We omit the details. \ \ \ \ \fbox{}

\begin{appendix}
\def\theequation{\Alph{section}.\arabic{equation}}

\section{Appendix}
\setcounter{equation}{0}

\subsection{Some basic results of regular variation theory}

In this subsection, we gather some basic results of regular variation regular variation theory that are needed in this paper. In most cases, we refer the reader to the basic references (such as \cite{BGT87}) and omit the proofs. However, in certain instances, we feel that we need to provide the proofs as they are not readily available in the literature.

\begin{prop}\lbl{p2.4} {\rm (}Representation Theorem{\rm)} ~The function $L(u)$
is slowly varying at infinity if and only if it can be written as
  $$L(u)=M(u)\exp\left\{\int_b^u\frac{\varphi(t)}t{\rm d}t\right\},
  \ \ \ \forall \ u\geq b$$
for some $b>0$, where the function $\varphi\in C([b,\infty))$ satisfies
$\dd\lim\limits_{u\to\infty}\varphi(u)=0$ and $M(u)$ is measurable on
$[b,\infty)$ with $\dd\lim\limits_{u\to\infty}M(u)=M^*\in(0,\infty)$.
\end{prop}

\begin{defi}\lbl{d2.3} A function ${\hat L}(u)$ is referred to as
{\em normalized slowly varying} at infinity if it satisfies the
requirements in Proposition {\rm \ref{p2.4}} with $M(u)$ replaced by
$M^*$. The function $R(u)=M^*u^\rho {\hat L}(u)$ is called {\em
normalized regularly varying} at infinity of index $\rho$
and we write $R(u)\in NRV_\rho$.

We say that $R(u)$ is regularly varying at the origin {\rm(}from the
right{\rm )} of index $\rho\in\mathbb R$, denoted by $R\in
RVZ_\rho$, if $R(1/u)\in RV_{-\rho}$. The set of all normalized
regularly varying functions at the origin of index $\rho$ is denoted
by $NRVZ_\rho$.
\end{defi}

Similar to \cite[Remark 2.2]{CD05}, we have

\begin{lem} \label{r2.3a}
Let $k\in \mathcal {K}_\ell$ with $\ell>0$. Then $k(1/u)$ belongs to
$NRV_{(\ell-1)/\ell}$. Furthermore, we have $K(1/u)$ belongs to
$NRV_{-1/\ell}$. And hence, $K(u)\in NRVZ_{1/\ell}$ and $ k(u)\in NRVZ_{(1-\ell)/\ell}$.
\end{lem}

\begin{lem} \label{l2.2} {\rm(\cite[Lemma 2.3]{CW11})} ~Suppose that the condition $(F_1)$ holds and the function
$\phi$ is given  by {\rm(\ref{1.10})}.  Then

{\rm (i)}~$-\phi'(t)=(p'F(\phi(t)))^{1/p}$, and
$|\phi'(t)|^{p-2}\phi''(t)=\frac{p'}{p}f(\phi(t))$,
where $p'=\frac{p}{p-1}$;

{\rm (ii)}~ $-\phi'\in NRVZ_{-r}$, $\phi\in NRVZ_{1-r}$,  where
$r=(\rho+1)/(\rho+1-p)>1$.
\end{lem}

Following the discussion of \cite[Section 2]{CD05},  we can prove

\begin{lem}\label{l2.3} \ A function $f\in NRV_\rho$ $($or $f \in NRVZ_\rho)$ if
and only if $f\in C^1[a_1, \infty) \ ($or $f\in C^1(0, a_1))$ for
some $a_1>0$ and $\lim\limits_{s \rightarrow \infty}\frac{sf'(s)}{f(s)}=\rho$ $($or $\lim\limits_{s \rightarrow 0^+}\frac{sf'(s)}{f(s)}=\rho)$.
\end{lem}

In view of Lemma \ref{l2.3} we can prove

\begin{lem}\label{l2.4} \ Assume that $f\in NRV_\rho \ ($or $f \in NRVZ_\rho)$ and $\rho\not =0$. If $f'(s)>0$, then the inverse function $f^{-1}(y)$ of $f(s)$ belongs to $ NRV_{1/\rho}\ ($or $f^{-1}(y)\in NRVZ_{1/\rho})$. If $f'(s)<0$, then the inverse function $f^{-1}(y)$ of $f(s)$ belongs to $NRVZ_{1/\rho}\ ($or $f^{-1}(y)\in NRV_{1/\rho})$.
 \end{lem}

Since $\phi\in NRVZ_{1-r}$ and $r>1$, $K\in RVZ_{1/\ell}$, by Lemma \ref{l2.4} we have
$\phi^{-1}\in RV_{1/(1-r)},\, K^{-1}\in RVZ_{\ell}$. Thanks to $k\in RVZ_{(1-\ell)/\ell}$, it follows that
 \bes
 &k\circ K^{-1}\circ\phi^{-1}(s)\in RV_{\frac{1-\ell}{1-r}},&\lbl{b.0}\\
 &f^*(s)=\lk(k\circ K^{-1}\circ\phi^{-1}(s)\rr)^pf(s)\in RV_{\frac{p(1-\ell)}{1-r}+\rho}=RV_{q}.&\lbl{b.1}\ees

In the following, $C$ represents a generic positive constant which can differ from line to line.

\begin{lem}\label{l2.5} Let $\varrho>0$ be a constant. Assume that $f$ is a positive continuous function in $(0, \infty)$, and $f\in RV_{\gamma}$ with $\gamma>\max\{1,\,\varrho\}$. Let $w(x,t)$ be a positive function with positive lower bound $w_0$. Then for any given constant $C\geq 1$, there is a constant $\Lambda>0$, which depends on $w_0$, $C$ and $\gamma$, such that, for all $(x,t)$,
  \[C(\Lambda+\Lambda^\varrho)f(w(x,t))<f(\Lambda w(x,t)).\]
 \end{lem}

{\bf Proof}. \ Denote $\alpha=\max\{1,\,\varrho\}$ and choose $\sigma>0$ satisfying  $\gamma-\sigma>\alpha$. Choose $\Lambda_0\geq 2$ and $0<\varepsilon\ll 1$ with  $(1-\varepsilon)\Lambda_0^\sigma>2C$. By the definition, there is a constant $z_0=z_0(\Lambda_0)>0$ such that\bes
 f(\Lambda_0z)\geq(1-\varepsilon)\Lambda_0^\gamma f(z), \ \ \ \forall \ z\geq z_0.
 \lbl{b.2}\ees
For any positive integer $j\geq 2$ and $z\geq z_0$, we have $\Lambda_0^{j-i}z\geq z\geq z_0$ for all $1\leq i\leq j-1$, and by inductively
 \bes
 f(\Lambda_0^jz)&=&f(\Lambda_0\Lambda_0^{j-1}z)\geq (1-\varepsilon)\Lambda_0^\gamma f(\Lambda_0^{j-1}z)\geq\cdots\nonumber\\
 &\geq& [(1-\varepsilon)\Lambda_0^\gamma]^jf(z)=[(1-\varepsilon)\Lambda_0^\sigma]^jf(z) \Lambda_0^{(\gamma-\sigma)j}\nonumber\\
 &>&2C\Lambda_0^{(\gamma-\sigma)j}f(z).
 \lbl{b.3} \ees
As $f(z)\to\infty$ as $z\to\infty$, we can choose $z_0$ so large that $f(z)\leq f(z_0)$ for all $w_0\leq z\leq z_0$.

Since $f\in RV_\gamma$, we have $\dd\lim_{s\to\infty}s^{-\gamma+\sigma}f(s)=\infty$. For the given constant $A> 2Cf(z_0)/w_0^{\gamma-\sigma}$, there is a constant $S_0>1$ such that $f(s)\geq As^{\gamma-\sigma}$ for all $s\geq S_0$.
It is obvious that there is a constant $\Lambda^*>1$ such that $\Lambda z>S_0$ for all $\Lambda\geq \Lambda^*$ and $w_0\leq z\leq z_0$. Therefore,
 \bes
 f(\Lambda z)\geq A(\Lambda z)^{\gamma-\sigma}\geq Aw_0^{\gamma-\sigma}\Lambda^{\gamma-\sigma}
 \frac{f(z)}{f(z_0)}\geq 2C\Lambda^{\gamma-\sigma}f(z), \ \ \ \forall \
 \Lambda\geq \Lambda^*, \ \, w_0\leq z\leq z_0.\qquad
 \lbl{b.4}\ees

Note that $\Lambda_0\geq 2$, we can choose an integer $j\geq 2$ such that $\Lambda_0^j\geq \Lambda^*$. Take $\Lambda=\Lambda_0^j$, then our conclusion is true.
In fact, for those $(x,t)$ with $w(x,t)\leq z_0$, by (\ref{b.4})
 \[f(\Lambda w(x,t))>2C\Lambda^{\gamma-\sigma}f(w(z,t))>2C\Lambda^\alpha f(w(x,t)).\]
For those $(x,t)$ with $w(x,t)>z_0$, by (\ref{b.3}),
 \bess
 f(\Lambda w(x,t))&=&f(\Lambda_0^jw(x,t))>2C\Lambda_0^{(\gamma-\sigma)j}f(w(x,t))\\
 &=&2C\Lambda^{\gamma-\sigma} f(w(x,t))>2C\Lambda^\alpha f(w(x,t)).
 \eess
The proof is complete.\ \ \ \fbox{}

\vskip 8pt When $f$ is increasing, the following better result can be obtained:

\begin{lem}\label{l2.6}
Under the conditions of Lemma $\ref{l2.5}$, we further assume that $f$ is increasing in $(0, \infty)$. Then for any given constant $C\geq 1$, there is a constant $\Lambda^*>0$, which depends on $w_0$, $C$ and $\gamma$, such that for all $(x,t)$ and all $\Lambda\geq \Lambda^*$,
  \[C(\Lambda+\Lambda^\varrho)f(w(x,t))<f(\Lambda w(x,t)).\]
\end{lem}

\begin{lem}\label{l2.7}
Suppose that $f\in RV_{\gamma}$ with $\gamma\in \mathbb{R}$, is
continuous, increasing and positive in $(0,\infty)$. Then for any given $\tau>0$, there is a constant $c=c(\tau)>0$, such that
  \[f(a)+f(b)>cf(a+b),\ \ \forall \ a,b\geq\tau.\]
\end{lem}

{\bf Proof}. ~We assume by contradiction that there exist two
sequences $\{a_n\}$ and  $\{b_n\}$, such that
$f(a_n)+f(b_n)\leq \frac{1}{n}f(a_n+b_n)$.
Since $f$ is continuous, increasing and positive in $[\tau,\infty)$, it follows
$a_n+b_n\to \infty$ as $n\to \infty$. Suppose that $a_n\leq b_n$
without loss of generality. Since $f$ is increasing on $(0,\infty)$,
we have $f(b_n)\leq \frac{1}{n}f(2b_n)$, i.e.,
$f(2b_n)/f(b_n)\geq n$.
On the other hand, since $f\in RV_{\gamma}$, we have
$\dd\lim_{n\to \infty}[f(2b_n)/f(b_n)]=2^\gamma$,
which is a contradiction.\ \ \ \fbox{}

\begin{lem}\label{l2.8}
Assume that $f$ is a positive continuous function and $f\in RV_{\gamma}$ with $\gamma>0$. Let $\tau>0$ be a given constant. Then there exist a positive, continuous and increasing function $g\in RV_{\gamma}$
and a constant $\sigma$, such that
  \[\sigma g(u)\leq f(u)\leq g(u), \ \ \ \forall \ u\geq\tau.\]
\end{lem}

{\bf Proof}. ~As $u^{-\gamma}f(u)$ is a slow variation function, by Proposition \ref{p2.4},
  $$u^{-\gamma}f(u)=M(u)\psi(u), \ \ \ \mbox{with} \ \psi(u)=\exp\left\{\int_b^u\frac{\varphi(t)}t{\rm d}t\right\}, \ \ \ u\geq b>0.$$
Since $\dd\lim\limits_{u\to\infty}M(u)=M^*\in(0,\infty)$, there is a constant $u_1>0$ such that $M^*/2<M(u)<2M^*$ for all $u\geq u_1$. Hence,
 \[\frac 12M^*u^\gamma\psi(u)\leq f(u)\leq 2M^*u^\gamma\psi(u),
  \ \ \ \forall \ u\geq u_1.\]
The direct computation gives $\lk(u^\gamma\psi(u)\rr)'=u^{\gamma-1}[\gamma+\varphi(u)]\psi(u)$.
By use of $\dd\lim\limits_{u\to\infty}\varphi(u)=0$, it follows that there is a $u_2>0$ such that $\gamma+\varphi(u)>0$ when $u\geq u_2$. That is, the function
$u^\gamma\psi(u)$ is increasing for $u\geq u_2$. Take a positive, continuous and increasing function $g_1(u)$ such that $g_1(u)=u^\gamma\psi(u)$ when $u\geq u_0=\max\{u_1, u_2\}$. It is obvious that $g_1(u)\in RV_\gamma$ and
$\frac 12 M^* g_1(u)\leq f(u)\leq 2M^*g_1(u)$ for all $u\geq u_0$.
Note that both $f$ and $g_1$ are continuous and positive in $[\tau, u_0]$, there is a constant $C>0$ such that
$C^{-1}g_1(u)\leq f(u)\leq Cg_1(u)$ for all $\tau\leq u\leq u_0$.
Take $g(u)=(2M^*+C)g_1(u)$, then
our conclusion holds. \ \ \ \fbox{}

\subsection{Some results on the unique solution of (\ref{1.5})}

\begin{lem}\lbl{l2.9} \ Assume that $g(u)$ and $h(u)$ are continuous functions  and $h(u)$ is  positive in $[a, \infty)$ for some constant $a>0$, and that $h(u)\in NRV_\gamma$ with $\gamma>1$. Let $v(t)$ and $w(t)$ be the positive solutions of the problems, respectively:
 \bess
  &v'(t)=- g(v(t)), \ \ t>0; \qquad v(0)=\infty,&\\
 &w'(t)=- h(w(t)), \ \ t>0; \qquad w(0)=\infty.&\eess
If $\lim\limits_{u\to\infty}\frac{g(u)}{h(u)}=1$, then
$\lim\limits_{t\to 0^+}\frac{v(t)}{w(t)}=1$.
 \end{lem}

{\bf Proof}. Note that $h(u)\in NRV_\gamma$ and $ \gamma>1+\nu$, similar to the proof of Lemma \ref{l2.8} we have that $h(u)/u^{1+\nu}$ is increasing when $u\geq u_0$ for some large constant $u_0$. Choose $\varepsilon_n$ with $0<\varepsilon_n<1/2$ and $\varepsilon_n\to 0$ as $n\to\infty$. In view of $g(u)/h(u)\to 1$ as $u\to\infty$, there is $u_n\geq 2^{1/\nu}u_0$ such that $(1-\varepsilon_n)h(u)\leq g(u)\leq(1+\varepsilon_n)h(u)$ when $u\geq u_n$.

Thanks to the fact that $v(t), w(t)\to\infty$ as $t\to 0$, there is $t_n>0$ such that $v(t)\geq u_n$ for all $0<t\leq t_n$. Therefore,
 \[(1-\varepsilon_n)h(v(t))\leq g(v(t))\leq(1+\varepsilon_n)h(v(t)), \ \ \forall \ 0<t\leq t_n.\]
Hence
 \[v'(t)\leq-(1-\varepsilon_n)h(v(t)), \ \ v(t)\geq 2^{1/\nu}u_0, \ \ \forall \ 0<t\leq t_n; \ \ \ v(0)=0.\]
Denote $\sigma_n=(1-\varepsilon_n)^{1/\nu}$. In view of $\varepsilon_n<1/2$, we see that $\sigma_nv(t)\geq u_0$ for all $n$ and $0<t\leq t_n$. Notice that $h(u)/u^{1+\nu}$ is increasing for $u\geq u_0$ and $\sigma_n<1$, it is easily seen that
$\sigma_n(1-\varepsilon_n)h(v(t))=\sigma_n^{1+\nu}h(v(t))\geq h(\sigma_nv(t))$ in $ (0, t_n]$. Therefore, $y_n(t)=\sigma_nv(t)$ satisfies
 \[y'_n(t)\leq-\sigma_n(1-\varepsilon_n)h(v(t))\leq-h(y_n(t)), \ \ \forall \ 0<t\leq t_n.\]
The comparison principle gives that $y_n(t)=\sigma_nv(t)\leq w(t)$ for all $0<t\leq t_n$. Hence $\lim\limits_{t\to 0^+}\frac {v(t)}{w(t)}\leq 1/\sigma_n$,
and consequently $\lim\limits_{t\to 0^+}\frac {v(t)}{w(t)}\leq 1$ by letting $n\to\infty$.

Similarly, we can prove that $\lim\limits_{t\to 0^+}\frac {v(t)}{w(t)}\geq 1$.\ \ \ \fbox{}

\vskip 4pt

The following corollary can be drawn; we shall omit the proof:

\begin{col}\lbl{c2.z1} \ In Lemma $\ref{l2.9}$, if we replace the condition $h(u)\in NRV_\gamma$  by $h(u)\in RV_\gamma$ with $\gamma>1$, the conclusion is also true.
\end{col}

Following the same line of argument, we can also establish:

\begin{lem}\lbl{l2.10} \ Assume that $g(u)$ and $h(u)$ are continuous functions and $h(u)$ is positive in $[a, \infty)$ for some constant $a>0$. Suppose further that $h(u)\in RV_\gamma$ with $\gamma>1$. Let $v(t)$ and $w(t)$ be the positive solutions of the problems, respectively:
 \bess
 &v'(t)=-g(v(t)), \ \ t>0; \qquad v(0)=\infty,&\\
 &w'(t)=-h(w(t)), \ \ t>0; \qquad w(0)=\infty.&
 \eess
If $\lim\limits_{u\to\infty}\frac{g(u)}{h(u)}=c$ for some constant $c>0$, then for any $\nu:0<\nu<\gamma-1$,
 \bess
 c^{1/\nu}\leq\liminf_{t\to 0^+}\frac{w(t)}{v(t)}\leq\limsup_{t\to 0^+}\frac{w(t)}{v(t)}\leq 1 \ \ \ \mbox{if} \ c\leq 1,\\
 1\leq\liminf_{t\to 0^+}\frac{w(t)}{v(t)}\leq\limsup_{t\to 0^+}\frac{w(t)}{v(t)}\leq c^{1/\nu}\ \ \ \mbox{if} \ c>1.
 \eess
 \end{lem}

Using Lemma \ref{l2.10}, we can also obtain:

\begin{lem}\label{l2.11} \ Assume that $g(u)$ and $h(u)$ are positive continuous differentiable functions, and $g(u)\in RV_\theta$, $h(u)\in RV_\gamma$ with $\gamma+\theta>1$. Let $c>0$ be a given constant. Denote by $v(t)$ and $w(t)$ solutions of the following problems, respectively:
 \bess
 &v'(t)=-g(c v(t))h(v(t)), \ \ t>0; \qquad v(0)=\infty,&\\
 &w'(t)=-g(w(t))h(w(t)), \ \ t>0; \qquad w(0)=\infty.&
 \eess
Then there exists a constant $C>1$ such that
$C^{-1}v(t)\leq w(t)\leq Cv(t)$ in $(0, T]$.
 \end{lem}

\subsection{Some results on the corresponding elliptic boundary blow-up problem}

In this final subsection, we recall, for the sake of ease of reference for the reader, some results about the
 boundary blow-up solutions of the $p-$Laplacian elliptic equation
  \xxx
  \Delta_p u= b(x)f(u),\ \ &x\in\Omega,\\
  u=\infty, \ \ &x\in\partial\Omega,
    \lbl{2.8}\zzz
where $b(x)\in C^\alpha(\Omega)$ for some
$0<\alpha<1$ with $b(x)\geq 0$ and $b(x)\not\equiv 0$ in $\Omega$. Set
$\Omega_0=\{x\in\Omega:\, b(x)=0\}$
and assume that ${\bar\Omega}_0\subset\Omega$ and $b(x)>0$ in
$\Omega\backslash\bar\Omega_0$.

\begin{theo}\lbl{t2.1}{\rm(\cite[Theorem 1.1]{CW11})} ~Assume that  $f$
satisfies  $(F_2)$ and
 \begin{quote}
 $(F_3)\ \ \displaystyle\int_1^\infty F^{-1/p}(t){\rm d}t<\infty$.
  \end{quote}
 Then the problem
{\rm (\ref{2.8})} has at least one positive solution.
\end{theo}

In fact, $(F_1)$ implies $(F_3)$, see \cite[Lemmas 2.1 and 2.2, and Remark 2.4]{CW11}.

\begin{theo} {\rm(\cite[Theorem 1.2]{CW11})}\lbl{t2.2} ~Assume that  $(F_1)-(F_2)$
hold, the function $\phi$ is defined by {\rm(\ref{1.10})}.

{\rm(i)} If there exist a function $k\in\mathcal{K}_\ell$ and a
positive continuous function $\beta(y)$ defined on $\partial\Omega$
such that
 \bess
   \lim_{\Omega\ni x\to y} \frac{b(x)}{k^p(d(x))}=\beta(y) \ \ {\rm uniformly \ for }
   \ y\in \partial\Omega,
  \eess
then {\rm (\ref{2.8})} has a unique positive solution and the blow-up
rate is given by
 \bess
 \lim_{\Omega\ni x\to y}\frac{u(x)}{\phi(K(d(x)))}
 =\left(\frac{r+\ell-1}{r\beta(y)}\right)^{\frac{r-1}{p}}
 \ \ \
 {\rm uniformly\ for}\ \ y\in\partial\Omega.
 \eess

{\rm(ii)} Suppose that there exist a function $k\in \mathcal {K}_\ell$
and constants $0<\beta_1\leq \beta_2$, $\delta>0$,
such that
 \bess
 \beta_1k^p(d(x))\leq b(x)\leq\beta_2k^p(d(x))\ \
   {\rm for \ all}\ \ x\in\Omega\ \ {\rm with}\ \ d(x)\leq \delta.
  \eess
Then the problem {\rm(\ref{2.8})} has a positive solution $u$
satisfying
 \bess
 \displaystyle\liminf_{d(x)\to 0}\frac{u(x)}{\phi( K(d(x)))}
 \geq \left(\frac{r+\ell-1}{r\beta_2}\right)^{\frac{r-1}{p}},\ \ \
 \displaystyle\limsup_{d(x)\to 0}\frac{u(x)}{\phi(K(d(x)))}
 \leq \left(\frac{r+\ell-1}{r\beta_1}\right)^{\frac{r-1}{p}}.
 \eess

{\rm(iii)} When $p=2$ and $\ell\not=0$, under the condition of {\rm(ii)},
the positive solution of {\rm(\ref{2.8})} is also unique.
\end{theo}

\begin{remark}\lbl{r2.2} \ From the proof of {\rm\cite[Theorem 1.2]{CW11}} it can be seen that, for any given $y\in \partial\Omega$, if the limit
  \bess
   \lim_{\Omega\ni x\to y} \frac{b(x)}{k^p(d(x))}=\beta(y)
  \eess
exists, then any positive solution $u(x)$ of {\rm(\ref{2.8})}
satisfies
 \[\lim_{\Omega\ni x\to y}\frac{u(x)}{\phi(K(d(x)))}
 =\left(\frac{r+\ell-1}{r\beta(y)}\right)^{\frac{r-1}{p}}. \]
\end{remark}

\end{appendix}

 \end{document}